
\mag 1200
\input amstex
\input amsppt.sty
\NoRunningHeads
\overfullrule 0pt

\let\al\alpha
\let\bt\beta
\let\gm\gamma 
\let\dl\delta \let\Dl\Delta
\let\epe\epsilon \let\eps\varepsilon \let\epsilon\eps

\let\la\lambda \let\La\Lambda
 \let\Om\Omega
 \let\phi\varphi

\def\O{\Cal O}
\def\o{\otimes}
\def\k{\frak k}

\def\d{\partial}
\def\<{\langle}
\def\>{\rangle}
\def\L{\Lambda}
\def\wo{\widetilde\o}

\def\C{\Bbb C}

\def\Z{\Bbb Z}

\def\bta{{\bar\eta}}

\def\bb{\frak b}

\def\F{\Cal F}

\def\g{\frak g}

  \def\Fwh{\,\wh{\!\F}}

\def\A{\roman A}

\def\E{\roman E}

\let\Fun\Fwh


\def\A{\Cal A} 

\def\Fun/{ {\text{Fun}}} 
\def\End/{ {\text{End}}} 
\def\Hom/{ {\text{Hom}}} 
\def\Rm/{\^{$R$-}matrix}

\def\M/{\Cal M } 
\def\A/{$A_{\tau,\eta}(sl_2)$} 
\def\Am/{$A_{\tau,\eta,\mu}(sl_2)$} 

\def\E/{$E_{\tau,\eta}(sl_2)$} 

\def\Ex/{ \widetilde {  E   }}
\def\EE/{ $\widetilde {  E   }_{\tau,\eta}(sl_2)$}

\def\eqg/{elliptic quantum group}
\def\Em/{$E_{\tau,\eta,\mu}(sl_2)$} 
\def\E/{$E_{\tau,\eta}(sl_2)$}

\def\Ce/{\C/{1\over \eta}\Z }

\def \g{\frak g}
\def \l{\frak l}
\def \h{\frak h}

\def \oh{\hat \otimes}
\def\End {{\text{End}}}
\def\Id  {{\text{Id}}}
\def\Alt  {{\text{Alt}}}
\def\coth  {{\text{cotanh}}}

\topmatter

\title
Geometry and classification of solutions
of the Classical Dynamical Yang-Baxter Equation
\endtitle

\author
Pavel Etingof and Alexander Varchenko
\endauthor

\thanks
The authors were supported in part by an 
NSF postdoctoral fellowship and
NSF grant
DMS-9501290
\endthanks

\date
June, 1997
\enddate

\abstract
{
The classical Yang-Baxter equation (CYBE) is an algebraic equation
central in the theory of integrable systems.
Its solutions were classified by Belavin and Drinfeld. Quantization
of CYBE led to the theory of quantum groups. A geometric interpretation
of CYBE was given by Drinfeld and gave rise to the theory of 
Poisson-Lie groups.

The classical dynamical Yang-Baxter equation (CDYBE) is an
important differential  equation analagous to CYBE and introduced by 
Felder as the consistency condition for the differential 
Knizhnik-Zamolodchikov-Bernard equations 
for correlation functions in conformal field theory on tori.
Quantization of CDYBE allowed Felder to introduce an interesting elliptic 
analog of quantum groups.  It becomes clear that numerous important
notions and results connected with CYBE have dynamical analogs.

In this paper we classify solutions to CDYBE and give geometric 
interpretation to CDYBE. The classification and interpretation are 
remarkably analogous to the Belavin-Drinfeld picture.
}
\endabstract

 \leftheadtext{ P.Etingof and A.Varchenko }
  \rightheadtext{ }
\endtopmatter

\document

\head 0.Introduction
\endhead

{\bf 0.1.} 
In 1984 Knizhnik and Zamolodchikov showed that correlation functions 
of conformal blocks 
on $\Bbb P^1$ for the Wess-Zumino-Witten (WZW) conformal field theory
for a simple Lie algebra $\g$ 
satisfy the differential equations
$$
\kappa \frac{\d F}{\d z_i}=\sum_{j\ne i}\frac{\Omega^{ij}}{z_i-z_j}F,
 \qquad i=1,...,N,\tag 1
$$
where $F$ is an analytic function of $N$ complex variables $z_1,...,z_N$
with values in $V_1\o...\o V_N$, $V_1,...,V_N$ are representations of $\g$, 
$\Omega^{ij}$ is the Casimir operator $\Omega\in (S^2\g)^\g$ acting in the
i-th and the j-th component, and $\kappa$ is a complex number. 
Equations (1) are called the Knizhnik-Zamolodchikov (KZ) equations. 
They play an important role in representation theory and mathematical physics.

Later Cherednik \cite{C} generalized the KZ equations: he considered
differential equations of the form
$$
\kappa \frac{\d F}{\d z_i}=
\sum_{j\ne i}r^{ij}(z_i-z_j)F,\tag 2
$$
where $r(z)$ is a meromorphic function on $\Bbb C$ with values in 
$\g\o \g$ satisfying the unitarity condition $r(z)=
-r^{21}(-z)$. Cherednik pointed out that system (2) is consistent iff
$r(z)$ satisfies the classical Yang-Baxter equation (CYBE): 
$$
[r^{12}(z_1-z_2),r^{13}(z_1-z_3)]+
[r^{12}(z_1-z_2),r^{23}(z_2-z_3)]+
[r^{13}(z_1-z_3),r^{23}(z_2-z_3)]=0.\tag 3
$$
In particular, for the simplest Yang's solution of (3), $r(z)=
\frac{\Omega}{z}$, we get the KZ equations. 

The geometric meaning of equation (3) was found by Drinfeld \cite{Dr}. 
Namely, he showed that any solution $r$ of (3), independent of $z$, 
and satisfying the condition $r+r^{21}\in (S^2\g)^\g$, defines 
a natural Poisson-Lie structure on a Lie group $G$ whose Lie algebra is 
$\g$. He also found that $z$-dependent solutions of (3) define analogous
structures on the loop group $LG$. 

If $\g$ is a simple Lie algebra, then equation (3) has many interesting 
solutions, both z-independent and z-dependent. These solutions, 
satisfying an additional nondegeneracy condition, were classified 
by Belavin and Drinfeld \cite{BD1},\cite{BD2}. In the z-dependent case, 
three types of solutions were found: rational, trigonometric, and elliptic 
(in terms of $z$).    

{\bf 0.2.} In \cite{B} Bernard found that 
correlation functions for WZW conformal blocks on an 
{\it elliptic curve} satisfy the differential equations
$$
\kappa \frac{\d F}{\d z_i}=
\sum_{j\ne i}r_{KZB}^{ij}(\la, z_i-z_j) F
+\sum_lx_l^i\frac{\d F}{\d x_l},\tag 4
$$
where $\la\in \h^*$, $\h$ is a Cartan subalgebra of $\g$, 
$F$ is an analytic function of complex variables $z_1,...,z_N$
and of a point $\la$, with values in $(V_1\o...\o V_N)^\h$,
$r_{KZB}(z,\la)$ is a 
particular meromorphic (in fact, elliptic) 
function with values in $(\g\o\g)^\h$, 
and $\{x_l\}$ is a basis of $\h$, which is also regarded as a linear 
system of coordinates on $\h^*$.   
Equations (4) are called the Knizhnik-Zamolodchikov-Bernard equations. 

Equations (4) are not of type (2), since they contain derivatives 
with respect to 
$x_l$ on the right hand side. Therefore, it is not surprising 
that $r_{KZB}$ does not satisfy the classical Yang-Baxter equation.
However, as was pointed out by Felder
\cite{F1}, $r_{KZB}$ satisfies a generalization of
the classical Yang-Baxter equation:
$$
\gather
\sum_lx_l^{(1)}\frac{\d r^{23}(z_2-z_3)}{\d x_l}
+\sum_lx_l^{(2)}\frac{\d r^{31}(z_1-z_3)}{\d x_l}+
\sum_lx_l^{(3)}\frac{\d r^{12}(z_1-z_2)}{\d x_l}+
\tag 5
\\
[r^{12}(z_1-z_2),r^{13}(z_1-z_3)]+
[r^{12}(z_1-z_2),r^{23}(z_2-z_3)]+
[r^{13}(z_1-z_3),r^{23}(z_2-z_3)]=0.
\endgather
$$
Moreover, this equation is a necessary and sufficient condition
for (4) to be consistent (for an arbitrary meromorphic function 
$r(\la, z)$ satisfying the unitarity condition $r(\la, z)=-r^{21}(\la, -z)$). 
This equation is called the {\it classical dynamical  
Yang-Baxter equation} (we will abbreviate it as 
CDYBE, or CDYB equation). The word ``dynamical'' refers to the fact 
that (5) is a differential rather than an algebraic equation, so it reminds 
of a dynamical system.  

This paper has two goals:

1. To exhibit the geometric meaning of CDYBE, similarly to 
Drinfeld's interprtation of CYBE
within the framework of the theory of Poisson-Lie groups. 

2. To classify solutions of CDYBE for simple Lie algebras and Kac-Moody 
algebras (over $\Bbb C$). 

{\bf 0.3.} The first goal is attained in Chapters 1 and 2. 
Namely, there we consider solutions of (5), 
for any pair of finite-dimensional Lie algebras $\h\subset \g$, 
which are independent of $z$, $\h$-invariant,  
and satisfy the generalized unitarity condition:
$r+r^{21}$ is a constant, invariant element of $S^2\g$.
To give geometric meaning to such solutions, we define and
study the notion of a {\it dynamical Poisson groupoid}
which is a special case of the notion of a Poisson groupoid
introduced by Weinstein \cite{W}.  
We show that any solution of (5) of the above type
naturally defines a dynamical Poisson groupoid,
which, as a manifold, equals $U\times G\times U$, where $U\subset \h^*$ 
is an open set.  
This construction illustrates how equation (5) arises naturally 
in the theory of Poisson groupoids. When $\h=0$, this construction 
reduces to the Drinfeld construction of a Poisson-Lie group 
from a solution of CYBE. 

As in the case of the usual CYBE, z-dependent solutions of (5) 
define analogous structures on the loop group $LG$. 

{\bf 0.4.} The second, more technically challenging goal, is 
(partially) attained in Chapters 3 and 4. 

In Chapter 3, we consider
z-independent solutions of CDYBE
which satisfy the condition that
$r+r^{21}$ is a constant invariant element, when $\g$ is a simple 
finite-dimensional Lie algebra, or, more generally, a Kac-Moody algebra. 
In this case, 
$r+r^{21}=\epe \Omega$, where $\Omega$ is the Casimir element, and
$\epe$ is a number called the coupling constant. 
For a simple Lie algebra $\g$, we classify 
all solutions. It turns out that there are two types of solutions 
-- rational (with zero coupling constant) and trigonometric
(with nonzero coupling constant). 

If $\g$ is an arbitrary Kac-Moody algebra, we classify 
solutions satisfying some additional conditions. Again, we find two types 
of solutions -- rational and trigonometric. 

We also classify invariant solutions of CDYBE
for a pair of Lie algebras ${\frak l}\subset\g$, where
$\g$ is a finite-dimensional simple Lie algebra, 
and ${\frak l}$ is a reductive Lie subalgebra of $\g$ containing 
the Cartan subalgebra $\h$. The classification is 
obtained by reduction of CDYBE for the pair 
${\frak l}\subset \g$ to CDYBE for the pair $\h\subset \g$. 
  
In Chapter 4 we are concerned with z-dependent solutions
of CDYBE. We consider such solutions, satisfying the unitarity condition and 
the condition 
that the residue of $r(\la, z)$ at $z=0$ equals $\epe\Omega$. 
As before, $\epe$ is called the coupling constant. 
We  classify all solutions with nonzero coupling constant.
It turns out that there are three types of solutions -- rational, 
trigonometric, and elliptic. 

We also explain that 
z-independent solutions for the affine Lie algebra $\tilde\g$
 can be interpreted as 
z-dependent solutions of
CDYBE for $\g$ with $\epe\ne 0$. 

{\bf 0.5.} The CDYB equation has a quantum analogue. 
This quantum equation is called the quantum 
dynamical Yang-Baxter equation (QDYBE). It was first introduced
by Gervais and Neveu \cite{GN} and later by Felder \cite{F1}, as a 
quantization of (5). This equation has important applications in the theory
of integrable systems \cite{ABB}.

QDYBE is a generalization 
of the usual quantum
Yang-Baxter equation, and it gives rise to the notion 
of a dynamical Hopf algebroid (or dynamical 
quantum groupoid), in the same way as 
the usual quantum Yang-Baxter equation gives rise to the notion of 
a Hopf algebra (quantum group). The notion of a dynamical Hopf algebroid
is a quantization of the notion of a dynamical Poisson groupoid, 
discussed in this paper. An example of a dynamical Hopf algebroid,
which is not a Hopf algebra, 
is the elliptic quantum group of \cite{F1,FV}.

In the following papers we will 
define and study the notion of a dynamical Hopf algebroid, 
give examples, and show how such objects arise naturally in 
representation theory of affine Lie algebras and quantum groups, and 
in conformal field theory. 

\subhead Acknowledgements\endsubhead

We thank Giovanni Felder who raised the question of geometric
interpretation of the classical dynamical Yang-Baxter equation. 
We are grateful to Olivier Schiffmann for careful reading of this manuscript
and pointing out a number of errors. We thank Vladimir Drinfeld 
and Jiang-Hua Lu for useful discussions. 

\vskip .1in 

\head 1. Dynamical Poisson groupoids.\endhead

\subhead 1.1. Definition of a dynamical Poisson groupoid\endsubhead

Recall that a groupoid is a
category where all morphisms are isomorphisms.
In this paper, we consider only groupoids
whose objects and morphisms form a set, 
and not only a class (i.e. groupoids built on small categories). 

Thus, a groupoid is defined by the following data: a set $X$ 
(of morphisms, or arrows) called the groupoid itself, 
a set $P$ (of objects) called the base of $X$, two surjective maps
$s,t:X\to P$ (sourse and target), a map
$m: \{(a,b)\in X\times X: s(a)=t(b)\}\to X$
(called the arrow composition map), and an injective map 
$E: P\to X$ (the identity morphism; $E(p)=id_p$), satisfying 
a number of obvious conditions. One of these conditions is 
the existence of an involution
$i:X\to X$ defined by the conditions
$s(i(x))=t(x)$, $s(x)=t(i(x))$, $m(x,i(x))=id_{t(x)}$, 
$m(i(x),x)=id_{s(x)}$. 
For brevity, when we talk about a particular groupoid, we will 
refer to its set of morphisms  $X$. 

A groupoid with one object is a group. Thus, the notion of a groupoid is
a generalization of the notion of a group. 

The role of the unit in a groupoid is played by 
the map $E$, and the role of inversion by $i$. 

A Lie groupoid is a groupoid with a smooth structure
(the sets of objects and morphisms are smooth manifolds, 
the structure maps are smooth, and some additional conditions,
see \cite{M}). 

According to \cite{W}, a
 Poisson groupoid is a Lie groupoid $X$ with a Poisson bracket such
that the graph of the composition map (defined 
only on a subset of $X\times X$) is a coisotropic submanifold of
$X\times X\times \bar X$, where $\bar X$ is the opposite Poisson manifold 
to $X$. For example, if $|P|=1$ (i.e. $X$ is a Lie group), the structure
of a Poisson groupoid is the same as a structure of a Poisson-Lie group. 

In this section we will define a class of Poisson groupoids which 
we call dymnamical Poisson groupoids. 

Let $G$ be a Lie group, and $\g$ its Lie algebra. 
Let $H$ be a connected Lie subgroup of $G$, and $\h$ the Lie algebra of $H$.
Let $U$ be an open subset of $\h^*$, which is invariant under the coadjoint 
action. 

Consider the manifold $X(G,H,U)=U\times G\times U$. 
This manifold has a natural structure of a Lie groupoid:
$X=X(G,H,U)$, $P=U$, $s(u_1,g,u_2)=u_2$, $t(u_1,g,u_2)=u_1$, 
$E(u)=(u,1,u)$, and $m((u_1,f,u_2),(u_2,g,u_3))=(u_1,fg,u_3)$. 
In the theory of groupoids, this groupoid is called 
the direct product of the trivial groupoid with base $U$ and 
the group $G$. 

The inversion of the groupoid $X$ is given by
$i(u_1,g,u_2)=(u_2,g^{-1},u_1)$. 

The manifold $X$ carries a pair of commuting actions of $H$:
the left action given by $l(h)(u_1,g,u_2)=(hu_1h^{-1},hg,u_2)$, and
the right action given by $r(h)(u_1,g,u_2)=(u_1,gh,h^{-1}u_2h)$
($huh^{-1}$ denotes, here and below, the coadjoint action of $h$ on $u$).  
The manifold $X\times X$ carries a (left) action of $H$
given by $\Delta(h)(x,y)=(r(h)^{-1}x,l(h)y)$. We will call 
this action the diagonal action. This action preserves the composition map. 

For $a\in\h$, let $a_1,a_2$ be the functions 
on $X$ given by $a_1(u_1,g,u_2)=a(u_1)$, $a_2(u_1,g,u_2)=
a(u_2)$.

Let $\{,\}$ be a Poisson bracket on $X$. 

\proclaim{Definition} The pair $(X,\{,\})$
 is said to be a dynamical Poisson groupoid if
the following two conditions are satisfied. 

(i) The actions $l,r$ are Hamiltonian, 
the maps $t,s$ are moment maps for them, and
for any $a,b\in \h$ one has $\{a_1,b_2\}=0$.

(ii) Let $X\bullet X:=X\times X//\Delta(H)$ be the Hamiltonian reduction of 
$X\times X$ by the diagonal action of $H$, and let
$\bar m: X\bullet X\to X$ be the reduction by $H$ 
of the composition map $m$ of $X$. Then 
$\bar m$ is a Poisson map. 
\endproclaim

{\bf Remark 2.} Let us explain condition (ii).
If condition (i) is satisfied, the diagonal action  
of $H$ is Hamiltonian, and $\mu(x,y)=t(y)-s(x)$ is a moment map for 
this action. Therefore, the set of zeros of the moment map, $\mu^{-1}(0)$, 
is precisely the domain of the map $m$. Thus, the definition of $\bar m$ makes 
sense. The space $X\bullet X$ equals
$U\times (Y/H)\times U$, where $Y=G\times U\times G$, and 
$H$ acts on $Y$ by $h\circ (f,u,g)=(fh^{-1},huh^{-1},hg)$.   
The space $X\bullet X$ has the Poisson structure of the Hamiltonian reduction. 
Therefore, it makes sense to require that the map $\bar m$ is Poisson. 

If $H=1$, then $\{,\}$ is a Poisson bracket on $G$. 
Condition (i) is vacuous, and condition (ii) says
that the multiplication map in $G$ is Poisson.
Thus, a dynamical Poisson groupoid with $H=1$ is simply 
a Poisson-Lie group. 

Let us compute the general form of the Poisson bracket on a dynamical Poisson 
groupoid $X=U\times G\times U$. 
Let $f$ be any function on $X$ which is pulled 
back from $G$. By the definition, we 
have the following Poisson 
commutation relations:
$$
\{a_1,b_2\}=0, \{a_1,b_1\}=-[a,b]_1,
\{a_2,b_2\}=[a,b]_2, \{a_1,f\}=R_af, \{a_2,f\}=L_af,\tag 1.1
$$
where $L_a,R_a$ are the left- and the right-invariant vector fields on
$G$ corresponding to $a$. So the only freedom that we have is in 
the bracket of functions on $G$. 

\subhead 1.2. The Hamiltonian unit\endsubhead

For Poisson Lie groups, it is known that the unit
$E:\{1\}\to G$ is a Poisson map. 
In the case of dynamical Poisson groupoids, this property fails: 
the image of $E$ is not a Poisson submanifold of $X$. 
Therefore it makes sense to extend $E$ so that its image
becomes the smallest Poisson manifold containing the image of $E$. 
This is done by introducing the Hamiltonian unit. 

Let $(T^*H)_U$ be the set of all $(h,p)\in T^*H$ ($h\in H,p\in T^*_hH$)
such that $h^{-1}p\in U$. We equip $(T^*H)_U$ with the standard 
symplectic structure.  

\proclaim{Definition} The Hamiltonian unit of a dynamical 
Poisson groupoid $X$ is the map $e: (T^*H)_U\to X$ given by
$e(h,p)=(ph^{-1},h,h^{-1}p)$. 
\endproclaim

It is easy to check that this map is Poisson, and 
its image is the smallest Poisson submanifold of $X$ containing 
the image of $E$. 

\subhead 1.3. Poisson groupoids and the CDYB equation\endsubhead

In this section we will construct examples of dynamical 
Poisson groupoids, and will be naturally led to the 
classical dynamical Yang-Baxter equation (CDYBE).

We work in the setting of Section 1.1.
Namely, we are considering the Lie groupoid $X=X(G,H,U)$. 
We want to make 
$X$ into a dynamical Poisson groupoid. 
As we know, the Poisson bracket on $X$ is partially defined by (1.1), and 
it remains to define 
the Poisson bracket of two functions pulled back from $G$. 

Recall \cite{Dr} that if $K$ is a Poisson Lie group, and $\k$ its Lie algebra, 
then a coboundary Poisson-Lie structure on $K$ is a Poisson-Lie structure 
with the Poisson bivector of the form 
$$
\pi=R(\rho)-L(\rho), \tag 1.2
$$
where $\rho\in \La^2\k$, and $L(\rho),R(\rho)$ are the left- and the
right-invariant tensor fields equal to $\rho$ at the identity. 
It is known \cite{Dr} that (1.2) defines a Poisson-Lie structure
iff 
$$
CYB(\rho):=[\rho^{12},\rho^{13}]+[\rho^{12},\rho^{23}]+[\rho^{13},\rho^{23}]
\in (\L^3\k)^\k.\tag 1.3
$$  

By analogy with this definition, we will look for a
Poisson bracket $\pi$ on $X$ such that for any functions $f_1,f_2$ 
pulled back from $G$
$$
\{f_1,f_2\}=(df_1\o df_2)(R(\rho(u_1))-L(\rho(u_2))),\tag 1.4
$$
where $\rho:U\to \L^2\g$ is a smooth function. 

For a Lie algebra $\g$ and a tensor $Z\in \g\o\g\o\g$, define 
$$
\Alt(Z)=Z^{123}+Z^{231}+Z^{312}.
$$ 

Let $\rho:U\to \L^2\g$ be a smooth function.
Then the differential $d\rho$ is a 1-form on $U$ with coefficients in
$\L^2\g$, so it can be considered as a function on $U$ with values in
$\h\o \La^2\g\subset \g\o\g\o\g$. 

Define the classical dynamical Yang-Baxter functional
$$
CDYB(\rho)=\Alt(d\rho)+CYB(\rho).\tag 1.5
$$

\proclaim{Theorem 1.1} Formulas (1.1),(1.4) define a Poisson structure 
on $X$ if and only if

(i) $\rho$ is $H$-equivariant, and

(ii) the element 
$Z=CDYB(\rho(u))$ is constant (in $u$) and lies in $(\L^3\g)^\g$.
\endproclaim

\demo{Proof} Property (i) is equivalent to the Jacobi identity 
for three functions $a_1,f_1,f_2$, where $a\in \h$, and $f_1,f_2$ 
are pulled back from $G$ (here, as usually, $a$ is regarded 
as a linear function on $U$). Also, when (i) is satisfied, 
then (ii) is equivalent to the Jacobi identity for three 
functions $f_1,f_2,f_3$ pulled back from $G$. 
\enddemo

Now let $\rho$ satisfy conditions (i),(ii).
Then $X$ equipped with the Poisson bracket $\{,\}$ defined by $\rho$ is  
a dynamical Poisson groupoid. Indeed, 
it is easy to see that the composition map $m: X\bullet X\to X$ 
given by $(u_1,g_1,u_2)\bullet (u_2,g_2,u_3)=
(u_1,g_1g_2,u_3)$ is Poisson. 

\proclaim{Definition} A dynamical Poisson groupoid defined by (1.4) 
is called a coboundary dynamical Poisson groupoid. 
\endproclaim

Recall \cite{Dr} that a coboundary Poisson-Lie group $K$ 
defined by (1.2) is called 
quasitriangular if it is equipped with a constant
element $T\in (S^2\g)^\g$, such that $CYB(\rho)=
\frac{1}{4}[T^{12},T^{23}]$. In this case 
the element $r=\rho+\frac{1}{2}T$ is a solution of the classical 
Yang-Baxter equation, and is called the classical r-matrix of $G$. 
If $T=0$, then the quasitriangular structure defined by $T$ is called
triangular. 

Thus, quasitriangular structures
on $G$ are parametrized by solutions $r$ of the CYB equation,
such that $r+r^{21}$ is a $\g$-invariant
element in $S^2\g$. Triangular structures are parametrized by
skew-symmetric solutions of the CYB equation. 

\proclaim{Definition} A coboundary dynamical Poisson groupoid is called 
quasitriangular if it is equipped with a  constant
element $T\in (S^2\g)^\g$, such that $CDYB(\rho)=
\frac{1}{4}[T^{12},T^{23}]$.
If $T=0$, the quasitriangular structure defined by $T$ is called
triangular.
\endproclaim

 In the quasitriangular case 
the function $r=\rho+\frac{1}{2}T$ is a solution of the classical 
dynamical Yang-Baxter equation
$$
CDYB(r)=0.\tag 1.6
$$
 and is called the classical 
dynamical r-matrix of $X$. Thus, quasitriangular structures
on $X$ are parametrized by solutions $r$ of the CDYB equation,
which are $\h$-invariant and such that $r+r^{21}$ is a constant $\g$-invariant
element in $S^2\g$. Triangular structures are parametrized by
skew-symmetric, $\h$-invariant solutions of the CDYB equation. 

{\bf Remark.} The material of this section trivially generalizes 
to the case when $G,H$ are complex Lie groups 
(algebraic groups, formal groups) rather than real Lie groups. 

\subhead 1.4. Gauge transformations of coboundary dynamical 
Poisson groupoids\endsubhead

Let $G$ be a complex Lie group, 
$H$ a commutative, connected complex Lie subgroup of $G$, 
$\g,\h$ their Lie algebras, and $U\subset \h^*$ a connected open set. 
Let $G^H$ be the centralizer of $H$ in $G$, and $\g^H$ its Lie algebra. 

In this section we will use the following notation: 
If $\alpha$ is a k-form on $U$ with values in a vector space $W$, 
then $\bar\alpha$ is the corresponding function 
$U\to\Lambda^k\h\o W$. 

Let $CB(G,H,U)$ be the set of all coboundary dynamical Poisson structures
on the groupoid $X=X(G,H,U)$. That is, $CB(G,H,U)$ is the set of all 
$\h$-invariant holomorphic functions $\rho: U\to \La^2\g$ such that
$CDYB(\rho)=Z\in (\Lambda^3\g)^\g$ is a constant. 

It turns out that there is a natural (infinite-dimensional) group 
which acts on $CB(G,H,U)$, such that 
the space of its orbits is finite-dimensional.  

Let $g:U\to G^H$ be a holomorphic function. Let $\eta=g^{-1}dg$ 
be a 1-form on $U$ with values in $\g^H$. The form 
$\eta$ defines a function $\bta: U\to \h\o\g^H$. 
For any function $\rho: U\to \Lambda^2\g$, define 
$$
\rho^g:=(g\o g)(\rho -\bta+\bta^{21})(g^{-1}\o g^{-1}).\tag 1.7
$$

\proclaim{Proposition 1.2} 
If $\rho\in CB(G,H,U)$ then $\rho^g\in CB(G,H,U)$, and $CDYB(\rho^g)
=CDYB(\rho)$.  
\endproclaim

\demo{Proof} Fix a basis $\{x_i\}$ of $\h$. We have 
$$
\gather
CDYB(\rho^g)=(\text{Ad}g)^{\o 3}\biggl(CDYB(\rho)+CDYB(\bta^{21}-\bta)\\
-[\rho^{12},\bta^{23}-\bta^{32}]-
[\rho^{12},\bta^{13}-\bta^{31}]-
[\rho^{13},\bta^{23}-\bta^{32}]+\\
[\rho^{23},\bta^{12}-\bta^{21}]+
[\rho^{13},\bta^{12}-\bta^{21}]+
[\rho^{23},\bta^{13}-\bta^{31}]+\\
Alt\biggl(\sum x_i\o [g^{-1}\frac{\d g}{\d x_i}\o 1+1\o 
g^{-1}\frac{\d g}{\d x_i},
\rho+\bta-\bta^{21}]\biggr)\biggr).\tag 1.8
\endgather
$$
Using the facts that $\rho$ is invariant under $\h$, and
 $CDYB(\rho)$ is invariant under $G$, we have 
$$
\gather
(Ad g^{-1})^{\o 3}(CDYB(\rho^g)-CDYB(\rho))=
CDYB(\bta^{21}-\bta)
+Alt([\rho^{23},\bta^{12}+\bta^{13}])+\\
Alt\biggl(\sum x_i\o [g^{-1}\frac{\d g}{\d x_i}\o 1+1\o 
g^{-1}\frac{\d g}{\d x_i},
\rho-\bta+\bta^{21}]\biggr).\tag 1.9
\endgather
$$
Simplifying the last two terms in (1.9), we get 
$$
\gather
(Ad g^{-1})^{\o 3}(CDYB(\rho^g)-CDYB(\rho))=\\
CDYB(\bta^{21}-\bta)+
Alt([\bta^{12}+\bta^{13},
\bta^{32}-\bta^{23}]).\tag 1.10\endgather
$$
However, since $\bta\in \h\o \g^H$, we have $[\bta^{12},\bta^{13}]=
[\bta^{12},\bta^{23}]=0$. Therefore, we have 
$$
\gather
(Ad g^{-1})^{\o 3}(CDYB(\rho^g)-CDYB(\rho))
=\\
CDYB(\bta^{21}-\bta)+
Alt([\bta^{12},
\bta^{32}]-[\bta^{13},\bta^{23}])=
Alt(d\bta^{21}-d\bta-[\bta^{13},\bta^{23}]).\tag 1.11
\endgather
$$
Let $\bar F_\eta:U\to \h\o\h\o\g^H$ be the function corresponding to 
the curvature form $F_\eta=d\eta+\frac{1}{2}[\eta,\eta]$. 
It is easy to see that the r.h.s. of (1.11) equals 
$-Alt(\bar F_\eta)$. Finally, observe that by the definition 
the form $\eta$ satisfies the zero-curvature condition 
$F_\eta=0$. 
Thus, $CDYB(\rho^g)=CDYB(\rho)$. 

$\square$\enddemo

It is easy to check that $(\rho^f)^g=\rho^{gf}$ for $f,g:U\to G^H$, 
so the assignment $\rho\to \rho^g$ defines a left action of the group
$\Sigma:=Map(U,G^H)$ of holomorphic functions on $U$ with values in $G^H$ 
on $CB(G,H,U)$. 

\proclaim{Definition} We will call transformations $\rho\to\rho^g$ 
the gauge transformations of the first kind. 
\endproclaim

{\bf Remark.} 
The gauge transformations of the first kind are 
especially simple if $g$ takes values 
in $H\subset G^H$. In this case $\rho^g=\rho+\bta^{21}-\bta$. 

Now let $\omega$ be a closed holomorphic 
2-form on $U$. This form defines a holomorphic function 
$\bar\omega: U\to\Lambda^2\g$. To this function there corresponds 
a transformation of $CB(G,H,U)$, given by $\rho\to \rho+\bar\omega$. 
We will call such transformations gauge transformations of the second kind. 

\proclaim{Proposition 1.3} If the form $\omega$ is exact, then the gauge 
transformation of the second kind $\rho\to\rho+
\bar\omega$ is also of the first 
kind. 
\endproclaim

\demo{Proof} Let $\xi$ be a 1-form on $U$ such that $\omega=d\xi$. 
This 1-form defines a function $\bar\xi: U\to \h$. Define a holomorphic 
function $g_\xi: U\to H$ by $g_\xi=e^{-\bar\xi}$. Then 
$\eta=g_\xi^{-1}dg_\xi=-d\bar\xi$. 
Thus, $\bta^{21}-\bta=\overline{d\xi}=\bar\omega$, as desired. 
$\square$\enddemo

 From now till the end of this section we will assume that $U$ is 
the formal polydisc. In this case by holomorphic functions 
we will mean arbitrary formal power series. Therefore, 
the constructions below make sense not only for the field $\C$ but for any 
field  of characteristic zero. 

Since $U$ is a formal polydisc, any closed form is exact. Thus, 
gauge transformations of the second kind are also of the first kind, and 
we will call them simply gauge transformations. 

Now we will show that the quotient space $CB(G,H,U)/\Sigma$ is 
finite-dimensional. 

\proclaim{Theorem 1.4} Let $\rho,r\in CB(G,H,U)$. Assume that 
the values of $\rho,r$ at the origin are equal, and 
$CDYB(\rho)=CDYB(r)$. Then $\rho=r^g$
for some $g\in\Sigma$. 
\endproclaim

\demo{Proof} 
Let $x_1,...,x_n$ be a basis of $\h$. We regard it as a system of
linear coordinates on $U$. The functions $r$, $\rho$ 
are by definition formal power series in $x_i$.

We will prove, by induction in $N$, that the statement 
of the theorem holds modulo 
terms of order $N+1$ of the power series. 
This is enough to prove the theorem. 

For $N=0$, the statement follows from the assumption of the theorem. 
Suppose we know it for $N=K$, and let us prove it for $N=K+1$.

We know that there exists a gauge 
transformation $g_K\in \Sigma$ such that the error term 
$\epe_K:=\rho-r^{g_K}$ is of order $K+1$. 
Let $\epe_K=E_K+\tilde\epe_K$, 
where $E_K$ is the part of $\epe_K$ of degree exactly 
$K+1$. 

Since $\rho$ and $r_K:=r^{g_K}$ satisfy the property 
$CDYB(\rho)=CDYB(r_K)$,
we have
$$
Alt(dE_K)=[CYB(r_K)-CYB(\rho)]_K,
$$
where $[*]_K$ denotes the degree $K$ homogeneous part of an expression $*$.
But according to our assumption, $\rho$ and $r_K$ coincide in 
degrees $\le K$, so we get
$$
Alt(dE_K)=0.
$$

Now we will study the equation $Alt(dE)=0$, which is the linearization 
of the CDYB equation near the zero solution. 

{\bf Lemma 1.} Let $E: U\to \La^2\g$ be a homogeneous polynomial function 
of degree $\ge 1$, 
invariant under $\h$, such that $Alt(dE)=0$. Then 
$E$ takes values in $\h\o\g^H+\g^H\o \h$.   

{\it Proof.} Because of the $\h$-invariance of $E$, 
it is enough to show that $E\in \h\o\g+\g\o\h$. 

Let $V\subset \g$ be 
a vector subspace which is 
a complement to $\h$. Since $E(u)\in \La^2\g$, we 
can write $E$ uniquely as a sum $E=E_{VV}+E_{V\h}+E_{\h V}+E_{\h\h}$, where
$E_{VV}\in \La^2V$, $E_{\h\h}\in \La^2\h$, $E_{V\h}\in V\o \h$, and
$E_{\h V}=-E_{V\h}^{21}$. 

The equation $Alt(dE)=0$ 
splits in 3 parts: the $\h VV$-part, the $\h\h V$-part, and the $\h\h\h$-part. 

The $\h VV$-part says that $dE_{VV}=0$, i.e. $E_{VV}$ is a constant. 
Since $E_{VV}$ is homogeneous of degree $\ge 1$, $E_{VV}=0$.  
The Lemma is proved. 

{\bf Lemma 2.} In the situation of Lemma 1, there exists a closed 1-form
$\eta$ on $U$ with values in $\g^H$, homogeneous of degree $K$,
 such that $E=\bta^{21}-\bta$. 

{\it Proof.} Let $V^H=V\cap \g^H$. Then $V^H$ is a complement to $\h$ in $V$. 
Let $\xi$ be the 1-form on $U$ with values in $V^H$ such that 
$\bar\xi=
-E_{\h V}$. 
The $\h\h V$-part of the equation $Alt(dE)=0$ says that
$d\xi=0$. 
Thus, $\xi$ satisfies the conditions of Lemma 2. 
Therefore, it is enough to prove Lemma 2 for $E'=E+\bar\xi^{21}-\bar\xi$. 

By the definition, $E'$ takes values in $\Lambda^2\h$, and $Alt(dE')=0$. 
Therefore, $E'$ defines a closed 2-form on $U$, which we denote by
$\omega$ (whe have $\bar\omega=E'$). Let $\zeta$ be a 1-form 
such that $\omega=d\zeta$, and let $\theta=d\bar\zeta$. 
Then $\bar\theta-\bar\theta^{21}=\bar\omega=E'$. The Lemma is proved. 

{\bf Lemma 3.} Let $\eta$ be a closed 1-form on $U$ with coefficients 
in a complex finite dimensional Lie algebra $\bb$, which has order $K$
at $0$. Then there exists a 1-form $\tau$ satisfying the 
zero-curvature equation 
$d\tau+\frac{1}{2}[\tau,\tau]=0$ such that $\tau-\eta$ is 
of order $\ge K+1$. 

{\it Proof.} Choose a function $\chi: U\to \bb$,
of order $K+1$, such that $\eta=d\chi$. 
Let $B$ be the Lie group of $\bb$, and consider the function
$g=e^\chi:U\to B$. Set $\tau=g^{-1}dg$. Then $\tau$ is the desired form.      

Now we return to the proof of the theorem. 
We start with the function $E=E_K$. Let $\eta$ be as in Lemma 2, 
and $\tau$ as in Lemma 3. Let $g:U\to G^H$ be the function such that 
$g^{-1}dg=\tau$. It is easy to see that for any 
$s: U\to \Lambda^2\g$ we have $s^g-s=E+\epe$, where $\epe$ 
are terms of order $K+1$ and higher. In particular, $r^{gg_K}-r^{g_K}$ 
equals $E$ modulo terms of order $\ge K+1$. Thus, if we set 
$g_{K+1}=gg_K$, we will get that $\rho-r^{g_{K+1}}$ is of order $\ge K+1$. 

The theorem is proved. 
$\square$\enddemo

\head 2. Biequivariant Poisson manifolds and groupoids.
\endhead

In this Chapter we will introduce the notion of an
$H$-biequivariant Poisson groupoid, which is a natural generalization of  
the notion of a dynamical Poisson groupoid. 

\subhead 2.1. Biequivariant Poisson manifolds\endsubhead

\proclaim{Definition} An $H$-biequivariant Poisson manifold
is a Poisson manifold $X$ equipped with two commuting, proper, free 
Hamiltonian actions of $H$ --
a left action $l:H\times X\to X$ and a right action $r: X\times H\to X$, 
and two maps $\mu_l,\mu_r: X\to \h^*$, which are moment maps for 
$l,r$, such that for any smooth
functions $a,b$ on $\h^*$ one has $\{a\circ \mu_l,b\circ \mu_r\}=0$.  
\endproclaim

{\bf Remark.} We remind (\cite{GHV},v.2, p. 135) that a smooth action 
of a group $H$ on a manifold $X$ is called proper if for any 
two compact sets $A,B\in X$ the set of elements of $a\in H$   
such that $B\cap aA\ne\emptyset$ is compact. It is known that 
if an action is proper and free then the quotient $X/H$ has 
a unique structure of a smooth manifold such that the natural map
$X\to X/H$ is a smooth submersion. 

We will denote the left and the right action of $h\in H$ on $x\in X$ by
$hx$ and $xh$, respectively. 

Observe that for any $H$-biequivariant 
Poisson manifold $X$ the maps $\mu_l,\mu_r$ are submersions, 
since the actions $l,r$ are free. 

Let $U\subset \h^*$ be an open set
invariant under the coadjoint action. We will say that 
an $H$-biequivariant 
Poisson manifold $X$ is over $U$ if the images of $\mu_l,\mu_r$ 
coincide with $U$. 

Let $\Cal C_U$ be the category of $H$-biequivariant Poisson manifolds 
over $U$, where morphisms are Poisson maps which commute
with $l,r,\mu_l,\mu_r$. We will now define the structure of a monoidal
category on $\Cal C_U$. 

Define the product $\bullet$ on $\Cal C_U$ as follows. 
Let $X_1,X_2\in \Cal C_U$, with 
actions of $H$ $l_1,r_1,l_2,r_2$, and 
moment maps $\mu_l^1,\mu_r^1,\mu_l^2,\mu_r^2$. 
Consider the left action $\Delta$ of $H$ on $X_1\times X_2$ 
by $\Delta(h)(x_1,x_2)=(x_1h^{-1},hx_2)$. The moment map for this action 
is $\mu_r^1-\mu_l^2$. Let $X_1\bullet X_2=X_1\times X_2//H$ be the Hamiltonian
reduction of $X_1\times X_2$ with respect to the action $\Delta$ of $H$. 
That is, $X_1\bullet X_2=Z/H$, where $Z\subset X_1\times X_2$ is the set
of points $(x_1,x_2)$ such that $\mu_r(x_1)=\mu_l(x_2)$.

It is easy to see that $Z$ is a smooth manifold 
(as $\mu_l,\mu_r$ are submersions), 
and  
$H$ acts freely and properly on $Z$, so $Z/H$ is smooth. 

The space $X_1\bullet X_2$ has a natural structure of an object of $\Cal C_U$:
it has a Poisson structure of Hamiltonian reduction, 
two free commuting actions of $H$ -- $l_1$ and $r_2$, and two 
moment maps $\mu_l^1$ and $\mu_r^2$ for them. Thus, $\bullet$ is a 
bifunctor $\Cal C_U\times \Cal C_U\to \Cal C_U$. 

It is easy to see that the operation $\bullet$ is associative:
$X_1\bullet (X_2\bullet X_3)=(X_1\bullet X_2)\bullet X_3=$\linebreak
$X_1\times X_2\times X_3//H\times H$, where $H\times H$ acts on
$X_1\times X_2\times X_3$ by $(h_1,h_2)(x_1,x_2,x_3)=
(x_1h_1^{-1},h_1x_2h_2^{-1},h_2x_3)$. 

Define an object $\bold 1\in \Cal C_U$ to be 
$(T^*H)_U$ (see Chapter 1), with the obvious left and right actions of $H$, and the moment maps
$\mu_l(p,h)=ph^{-1}$, $\mu_r(p,h)=h^{-1}p$, $h\in H$, $p\in T_h^*H$.
It turns out that $\bold 1$ is a unit object of $\Cal C_U$. 

Indeed, let us check that $\bold 1\bullet X$ is naturally isomorphic 
to $X$. 
We have $Z=\{(h,p,x):h^{-1}p=\mu_l(x)\}$. Thus, $Z$ is naturally 
identified with $H\times X$, via $(h,p,x)\to (h,x)$ 
The action of $H$ on $Z$ is by $h(h',x)=(h'h^{-1},hx)$. 
Thus, the quotient $Z/H$ is naturally isomorphic to $X$, and it is easy to 
check that the Poisson structure on $X$, the two actions of $H$, and the 
corresponding moment maps are the original ones. Similarly one
checks that $X\bullet \bold 1$ is naturally isomorphic to $X$. 
The unit object axioms are checked directly. Thus, 
$(\Cal C_U,\bullet,\bold 1)$ is a monoidal category with a 
unit object $\bold 1$.  

{\bf Remark.} Let $\Cal D$ be the category whose objects are 
manifolds with two commuting, free, proper actions of $H$ -- a left 
action $l$ and a right action $r$. This category has a natural 
monoidal structure, with product being the fiber product $\times_H$   
over $H$, and the unit object $H$ (with obvious $l$ and $r$). 
Then we have a natural functor $T^*$ from $\Cal D$ to $\Cal C_{\h^*}$
-- the functor of cotangent bundle: $M\to T^*M$. (Indeed, 
it is well known that if $H$ acts on $M$, then its induced 
action on $T^*M$ is Hamiltonian, with moment map
$\mu(m,p)(L)=\<L(m),p\>$, where $m\in M$, $p\in T_mM$, 
$L\in \h$, and $L(a),a\in M$ is the corresponding vector field on $M$.
Thus, $T^*M$ is naturally an object of $\Cal C_{\h^*}$.)     
The main property of the functor $T^*$ is that it is a monoidal functor: 
$T^*(M_1\times_H M_2)=T^*M_1\bullet T^*M_2$. 

Let $X\in \Cal C_U$. Denote by $\bar X$ the new object of $\Cal C_U$
obtained as follows: $\bar X$ is $X$ as a manifold, with the opposite Poisson 
structure $-\{,\}$, the left and the right actions of $H$ permuted
(i.e. the left, respectively right, action of $h$ on $\bar X$ is the right, 
respectively left, action of $h^{-1}$ on $X$), and the moment maps
also permuted. We will call $\bar X$ the dual object to $X$. 
By a reflection on $X$ we will mean an involutive morphism 
$i:X\to \bar X$. We will often write $x^{-1}$ for $i(x)$. 

Let $X\in \Cal C_U$ and $i:X\to\bar X$ be a reflection map.
Let $\phi^i_+,\phi^i_-: X\to X\times X$ be given by the formulas $\phi^i_+(x)= 
(x,x^{-1})$, $\phi^i_-(x)=(x^{-1},x)$. It is easy to see that these 
maps descend to maps (not necessarily Poisson) $\psi^i_\pm: X\to X\bullet X$, 
as $\mu_r(x^{-1})=\mu_l(x)$, $\mu_l(x^{-1})=\mu_r(x)$. 

\subhead 2.2. Biequivariant Poisson groupoids\endsubhead

\proclaim{Definition} Let $X\in \Cal C_U$. We will say that $X$ 
is an $H$-biequivariant Poisson semigroupoid if it is equipped 
with an associative morphism $m: X\bullet X\to X$
(called multiplication). In this case, a unit in $X$ is 
a morphism $e:\bold 1\to X$ such that the morphisms 
$m\circ (e\bullet id)$, $m\circ (id\bullet e)$ are the canonical morphisms
$\bold 1\times X\to X$, $X\times \bold 1\to X$. 
Further, such $X$ is called an $H$-biequivariant Poisson groupoid
if it is equipped with a reflection $i:X\to\bar X$ (called the inversion map), 
 such that $m(\psi^i_+(x))=e(1,\mu_l(x))$, 
$m(\psi^i_-(x))=e(1,\mu_r(x))$, where $(1,\mu_{l,r}(x))\in (T_1^*H)_U$.   
\endproclaim

{\bf Remark 1.} If $X$ has a unit, it is unique, as for any two units
$e_1,e_2$, we have $e_1=m\circ (e_1\bullet e_2)=e_2$. If the inversion map
$i$ on $X$ exists, it is also unique, as for any two inversion maps
$i_1,i_2$, $m_3(i_1(x)\bullet x\bullet i_2(x))=i_1(x)=i_2(x)$, where $m_3: 
X\bullet X\bullet X\to X$ is the multiplication map. 

{\bf Remark 2.} If $H$ is trivial, the notion of an $H$-biequivariant Poisson
groupoid coincides with the notion of a Poisson-Lie group. 

If $X$ is a dynamical Poisson groupoid, it is automatically an 
$H$-biequivariant Poisson groupoid. Indeed, in this case 
 $X$ is an $H$-biequivariant Poisson manifold,
with $\mu_l=t,\mu_r=s$. The composition map is the map 
$\bar m$ defined in Section 1.1, which is obviously associative. 
The unit axiom is satisfied for the Hamiltonian unit $e: \bold 1\to X$. 
Finally, it is easy to check that 
the inversion map is anti-Poisson and 
satisfies the inversion axiom. 

As we mentioned, 
the notion of a Poisson groupoid is known in the literature \cite{W}.
So let us justify the usage of this term in our paper, by showing that
$X$ is indeed a Poisson groupoid in the sense of \cite{W} in a natural way. 

Let $X$ be an $H$-biequivariant Poisson groupoid
over $U$. Then $X$ has a 
natural structure of groupoid in the usual sense. Namely, 
$P=U$, $s=\mu_r,t=\mu_l$, the composition map is
given by the multiplication map $m$, and $E=e|_{U}$, where 
$U=\{(1,p)\in (T^*H)_U\}$. It is easy to check that this groupoid 
is in fact a Lie groupoid in the sense of \cite{M}. 

\proclaim{Proposition 2.1} $X$ is a Poisson groupoid in the sense 
of \cite{W}. 
\endproclaim

\demo{Proof} We should show that the graph of composition is coisotropic. 
This follows from the following easy lemma.

{\bf Lemma.} Let $X$ be a Poisson manifold with a proper, 
free hamiltonian
action of a connected Lie group $H$ with moment map $\mu: X\to \h^*$. 
Let $X_0=\mu^{-1}(0)$. 
Let $Y$ be another Poisson manifold, and $f:X_0\to Y$ be an $H$-invariant 
smooth map. Then: $f$ descends to a Poisson map $X//H\to Y$ (where
$X//H$ is the hamiltonian reduction) if and only if the graph of $f$
is coisotropic in $X\times \bar Y$ (where $\bar Y$ is $Y$ with the opposite
Poisson structure).

{\it Proof.} Straightforward. 

To prove Proposition 2.1, it is enough to apply this Lemma to 
$f=m$, where $m$ is the multiplication map in the groupoid. 
\enddemo

{\bf Remark.} We have defined dynamical and
$H$-biequivariant Poisson 
groupoids in the category of smooth manifolds. 
Similarly one can define the same objects in the 
categories of complex analytic, 
formal, and algebraic varieties. In the formal and algebraic settings,
we can work over an arbitrary field of characteristic zero. 
These generalizations are straightforward, and we will not give them here. 

\subhead 
2.3. $H$-biequivariant Poisson algebras\endsubhead

In this and the next 
section we will sketch the constructions of the previous two sections 
in the algebraic language, i.e. working with Poisson algebras rather
than Poisson manifolds. This is related to  
the previous two sections by the operation of taking spectrum.

Let $k$ be a field of characteristic zero.
Let $A$ be a Poisson algebra
over $k$, $H$ a connected affine algebraic group, and
$\psi: H\times A\to A$ be a right algebraic
action of $H$ on $A$ by Poisson automorphisms.
(``Algebraic'' means that $A$ is a sum of finite dimensional 
representations of $H$). 

Let $\h$ be the Lie algebra of $H$.
Then the variety $\h^*$ has a natural Poisson structure. 
Let $U\subset \h^*$ be an $H$-invariant open set. A Poisson 
homomorphism $\mu: \O(U)\to A$ 
(where $\O(X)$ denotes the ring of algebraic functions on a variety $X$)
is called a {\it moment map} for $\psi$
if for any regular function $g$ on $U$, and any $f\in A$ we have 
$$
\{\mu(g),f\}=\sum_j\mu(\frac{\d g}{\d y_j})\cdot d\psi|_{h=1}(y_j,f).
$$ 
Here $y_j\in \h$ are a linear system of coordinates on $U$, $h\in H$, and
$d\psi|_{h=1}:\h\times A\to A$ is the differential of $\psi$ at $h=1\in H$.  
In particular, for a linear function
on $U$ given by $a\in \h$ the last equation is 
$$
\{\mu(a),f\}=\sum_j d\psi|_{h=1}(a,f).
$$ 
For a left action of $H$, a moment map is defined in the same way, 
with the only difference that it is anti-Poisson rather than Poisson. 

\proclaim{Definition} An $H$-biequivariant Poisson algebra over $U$
is a 5-tuple $(A,l,r,\mu_l,\mu_r)$, where 
$A$ is a Poisson 
algebra with $1$ over $k$, 
$l,r$ is
a pair of commuting algebraic 
 actions of $H$ on $A$ (a left action and a right action) 
by Poisson algebra automorphisms,
and $\mu_l,\mu_r: \O(U)\to A$ 
are moment maps for $l$, $r$, such that 

(i) $\mu_l,\mu_r$ are embeddings, and their images Poisson commute;

(ii) There exists an $l(H)\times r(H)$-invariant 
subspace $A_0^l$ of $A$ such that the multiplication map 
$\mu_r(\O(U))\o A_0^l\to A$ is a linear isomorphism; 
there exists an $l(H)\times r(H)$-invariant 
subspace $A_0^r$ of $A$ such that the multiplication map 
$\mu_l(\O(U))\o A_0^r\to A$ is a linear isomorphism. 

A morphism of $H$-biequivariant Poisson algebras over $U$ is a morphism
of Poisson algebras which preserves $l,r$ and $\mu_l,\mu_r$. 
\endproclaim

{\bf Remark 1.} Form $[l,r]=0$ it follows that 
$\{\mu_l\circ x,\mu_r\circ y\}$ is a central element (in the Poisson sense) 
for any $x,y\in \h$, but it does not follow that this commutator equals to 
zero. So we require that it is zero by condition (i). 

{\bf Remark 2.} Condition (ii) is of technical nature, and is not very 
important in the discussion below. 

Denote the category of $H$-biequivariant Poisson algebras over $U$ by 
$\Cal A_U$.

For convenience we will write $l(h)a$ as $ha$ and $r(h)a$ as $ah$. 

Let us now describe the monoidal structure on $\Cal A_U$.

Let $A,B\in \Cal A_U$. Let $l_A,r_A,l_B,r_B,
\mu_l^A,\mu_r^A,\mu_l^B,\mu_r^B$ be the corresponding 
actions and moment maps. Consider the action of the group $H$ in
$A\o B$ by $\Delta(h)(a\o b)=ah^{-1}\o hb$. We will construct a new 
$H$-biequivariant Poisson algebra $A\wo B$, which is obtained by 
Hamiltonian reduction of $A\o B$ by this action of $H$. 
 
Denote by $A*B$ the product $A\o_{\O(U)}B$,
where $\O(U)$ is mapped to $A$ via $\mu^A_r$ and to $B$ via $\mu^B_l$. 
The algebra $A*B$ has two commuting actions of $H$ ($l_A\o 1$ and $1\o r_B$). 
But we cannot claim that 
$A*B\in \Cal A_U$, since 
the Poisson structure on $A\o B$ does not, in general, descend to 
$A*B$. 

However, the action $\Delta$ of $H$ on $A\o B$ descends to
one on $A*B$, so we can define $A\wo B:=(A*B)^H$, where $H$ acts by $\Delta$. 
It is easy to check that the Poisson structure on $A\o B$ descends to
one on $A\wo B$ (Hamiltonian reduction). 
The two actions of $H$ and their moment maps also
descend to $A\wo B$. So, in order to check that $A\wo B\in \Cal A_U$, it 
suffices to check properties (i) and (ii). 

Using 
properties (i) and (ii) of the moment maps $\mu_l^A,\mu_r^A,\mu_l^B,\mu_r^B$,
it is easy to see that $A*B$ is naturally identified with
$\mu_l^A(\O(U))\o A_0^r\o B_0^r$, and $A\wo B$ is identified with
$\mu_l^A(\O(U))\o (A_0^r\o B_0^r)^H$, where $H$ acts by
$a\o b\to ah^{-1}\o hb$. This implies properties (i) and (ii) 
for the moment map $\mu_l^A\o 1: \O(U)\to A\wo B$, corresponding to
the left action of $H$ on $A\wo B$ (with $(A\wo B)_0^r=(A_0^r\o B_0^r)^H$). 
For the moment map $1\o \mu^B_r:\O(U)\to A\wo B$ corresponding 
to the right action, these properties are proved 
analogously. 

Thus, $A\wo B\in\Cal A_U$. 
It is clear that the assignment $A,B\to A\wo B$ is a bifunctor
$\Cal A_U\times \Cal A_U\to \Cal A_U$. 

Consider the algebra $\O((T^*H)_U)$, with the standard
Poisson structure, equipped with the standard actions $l,r$ 
of $H$ on left and right given by
$(x,p)\to (h_1xh_2,h_1ph_2)$. Let $M_l,M_r: (T^*H)_U\to U$
be given by $(h,p)\to ph^{-1}$, $(h,p)\to h^{-1}p$. Let  
$\mu_{l,r}=M_{l,r}^*:\O(U)\to \O((T^*H)_U)$. It is easy to check 
that $\mu_{l,r}$ are moment maps for $l,r$.  

Let $\bold 1=(\O((T^*H)_U),l,r,\mu_l,\mu_r)\in \Cal A_U$. 
It is easy to check that 
we have natural isomorphisms $A\wo \bold 1\equiv A\equiv \bold 1\wo A$.

\proclaim{Proposition 2.2} (i) $(A\wo B)\wo C=A\wo (B\wo C)$. 

(ii) $\bold 1$ is a unit object in $\Cal A_U$
with respect to $\wo$, and $(\Cal A_U,\wo,\bold 1)$ is
a monoidal category.
\endproclaim

\demo{Proof} Easy.
\enddemo

Let $A\in \Cal A_U$. Denote by $\bar A$ the new object of $\Cal A_U$
obtained as follows: $\bar A$ is $A$ as an algebra, with the opposite Poisson 
structure, the left and the right actions of $H$ permuted
(i.e. the left, respectively right, action of $h$ on $\bar A$ is the right, 
respectively left, action of $h^{-1}$ on $A$), and the moment maps
also permuted. We will call $\bar A$ the dual object to $A$. 
By a {\it reflection} on $A$ we will mean an involutive morphism 
$i:\bar A\to A$. 

Let $A\in \Cal A_U$ and $i:\bar A\to A$ be a reflection.
Let $\phi^i_+,\phi^i_-: A\o A\to A$ be given by the formulas $\phi^i_+(a\o b)= 
ai(b)$, $\phi^i_-(a\o b)=i(a)b$. It is easy to see that these 
maps descend to maps (not necessarily Poisson) $\psi^i_\pm: A\wo A\to A$. 

\subhead 2.4. $H$-biequivariant Poisson-Hopf algebroids\endsubhead

Now let us define the algebraic version of the notion of 
an $H$-biequivariant Poisson groupoid -- the notion of an
$H$-biequivariant Poisson-Hopf algebroid. 

\proclaim{Definition} Let $A$ be an $H$-biequivariant Poisson algebra.
Then $A$ is called an $H$-biequivariant Poisson-Hopf algebroid over $U$
if it is equipped with a coassociative $\Cal A_U$-morphism   
$\Delta: A\to A\wo A$ called the coproduct, 
an $\Cal A_U$-morphism $\epsilon: A\to \bold 1$ called the counit, and 
a reflection $S: \bar A\to A$ called the antipode,   
such that 

(i) $(id\bullet \epsilon)\circ \Delta=(\epsilon \bullet id)\circ \Delta
=id$, and 

(ii) $\psi_+^S\circ \Delta=\mu_l\circ P\circ \epsilon$,
$\psi_-^S\circ \Delta=\mu_r\circ P\circ \epsilon$, where $P: \bold 1\to
\O(U)$ acts by $f(x,p)\to f(1,p)$.  
\endproclaim 

{\bf Remark.} In the above discussion, $U$ is a Zariski open set. 
If $k=\Bbb R$ or $\Bbb C$, then we can take $U$ to be an open set 
in the usual sense, and define $\O(U)$ to be the algebra of smooth,
respectively analytic, functions on $U$. Then we can repeat sections 
2.3, 2.4 and thus define the notions of
an $H$-biequivariant Poisson algebra and Poisson-Hopf algebroid
over $U$. 
Similarly, one can take $U$ to be the infinitesimal neighborhood of zero 
in $\h^*$ (i.e. $\O(U)=k[[\h]]$). The material of Sections 2.3 and 2.4
can also be generalized to this case in a straightforward way.

The constructions of Section 1.3 
can easily be put in the algebraic framework of
Sections 2.3, 2.4. Let $G$ be an affine algebraic group, and 
$H$ a connected algebraic subgroup of $G$. 
Let $X(G,H,U)$ be a coboundary dynamical Poisson groupoid
defined by (1.1),(1.4). 
Consider the algebra $A=\O(U)\o \O(G)\o \O(U)$, where 
$\O(G)$ is the algebra of polynomial functions on $G$.
It is easy to see that $A$ is closed under the Poisson bracket
defined by (1.1),(1.4), 
so it is a Poisson algebra. 
The two actions of $H$ on $A$ are defined by
$$
\gather
l(h)[a\o f\o b](u_1,g,u_2)=[a\o f\o b](h^{-1}u_1h,h^{-1}g,u_2),\\
r(h)[a\o f\o b](u_1,g,u_2)=[a\o f\o b](u_1,gh^{-1},hu_2h^{-1}),\tag 2.1
\endgather
$$
and the corresponding moment maps are the maps corresponding to the 
projections of $X$ to the first and second component of $U$. 
The coproduct, counit, and antipode in $A$ are defined by the groupoid 
structure on $X$. Thus, $A$ is an $H$-biequivariant Poisson-Hopf algebroid.
We will call such a $H$-biequivariant Poisson-Hopf algebroid a dynamical
Poisson-Hopf algebroid. 

The notion of a dynamical Poisson-Hopf algebroid will be useful 
to us in the next paper.

\head 3. Classification of classical dynamical $r$-matrices
\endhead

\subhead 3.1.
Kac-Moody algebras [K, Ch. 2]
\endsubhead

Let $A = (a_{i,j} )_{i,j=1}^n$ be a symmetrizable generalized Cartan matrix.
Let $\g(A)$ be the associated Kac-Moody Lie algebra, $\h$ its Cartan
subalgebra,
$$
\g= \h \oplus \oplus _{\al \in \Delta} \g_\al 
$$
the root decomposition of the Kac-Moody Lie algebra. 
Let $(\cdot,\cdot)$ be an invariant nondegenerate bilinear
 form on $\g$ and $\Omega \in \g\oh \g$
the associated Casimir operator.
Here $\g\oh \g$ denotes the completed tensor product.
\subhead
Remark
\endsubhead
In this work we consider the Kac-Moody Lie algebras associated with a 
symmetrizable generalized Cartan matrix, although all theorems with the same proof
are valid for more general Kac-Moody Lie algebras associated with a symmetric
complex matrix, see [SV], [V, 11.1.10].

\subhead 3.2. Classical dynamical r-matrices
\endsubhead

Let $\g $ be a Kac-Moody Lie algebra, $\h$ its Cartan subalgebra.
A meromorphic function
$$
r \, : \, \h^* \, \to \g \oh \g
$$
is called {\it a classical dynamical r-matrix associated with
the pair $\h \subset \g$} if it satisfies the following
three conditions.
\item{1.}{{\it The zero weight condition,} 
$$
[h\otimes 1 + 1\otimes h\, , \, r(\lambda)]=0
\tag 3.1
$$
for any $\la \in \h^*$ and $h\in \h$.}
\item{2.}{{\it The generalized unitarity,} 
$$
r^{12}(\la) + r^{21}(\la) = \epe \, \Omega
\tag 3.2
$$
for some constant $\epe \in \C$ and all $\la$.}
\item{3.}{{\it The classical dynamical Yang-Baxter equation, CDYB,} 
$$
\Alt (d r) \, + \, [r^{12},r^{13}]\,+ \, [r^{12},r^{23}]\,
+\,[r^{13},r^{23}] \, = \, 0 \, .
\tag 3.3
$$
}

Here we use the following notations and conventions.
By a meromorphic function
$r \, : \, \h^* \, \to \g \oh \g$ satisfying the zero weight condition
 we mean a function of the form
$r\,=\,r_\h\, +\, \sum_{\al \in \Dl}\,r_\al$ where
$\Dl$ is the set of roots and
$r_\h \, : \, \h^* \, \to \h \otimes \h$, \,
$r_\al \, : \, \h^* \, \to \g_\al \otimes \g_{-\al}$
are meromorphic maps of finite dimensional spaces.

 If $X\in\End(V_i)$,
then
we denote by $X^{(i)}\in\End(V_1\otimes\dots\otimes V_n)$
the operator $\cdots\otimes\Id\otimes X\otimes\Id\otimes\cdots$, acting
non-trivially on the $i$th factor of a tensor product of vector spaces,
and
if $X=\sum X_k\otimes Y_k\in\End(V_i\otimes V_j)$, then we set
$X^{ij}=\sum X_k^{(i)}Y_k^{(j)}$.   

The differential of the r-matrix is considered in (3.3) as a 
meromorphic function
$$
dr\, : \, \h^* \, \to \g \hat\otimes\g \hat\otimes\g , \qquad
\lambda \, \mapsto \, \sum _i \, x_i \otimes {\partial r
\over \partial x_i}(\la)\, , 
$$
where $\{x_i\}$ is any basis in $\h$. We denote by $\Alt (dr)$ the following
symmetrization of $dr$ with respect to even permutations of 
1, 2, 3,
$$
\Alt (dr) \, = \, \sum _i \, x_i^{(1)} {\partial r^{23} 
\over \partial x_i}\, + \, \sum _i \, x_i^{(2)} {\partial r^{31} 
\over \partial x_i}\, + \, \sum _i \, x_i^{(3)} {\partial r^{12 }
\over \partial x_i}\, .
$$
The CDYB equation is an equation in $\g \hat\otimes\g \hat\otimes\g$.

The constant $\epe$ in (3.2) is called {\it the coupling constant}.


\subhead 3.3. Classification of the classical dynamical r-matrices with
nonzero coupling constant
\endsubhead

\proclaim {Theorem 3.1}

\item{1.}{ Let $x_i, \, i=1,...,N,$ be a basis in $\h$. For
any positive root $\al \in \Delta_+$ let
$e^i_\al, \, i=1,...,N_\al,$ be a basis in the root space $\g_\al$
and $e^i_{-\al}, \, i=1,...,N_\al,$  the dual 
basis in $\g_{-\al}$. Let $\nu$ be
an element in $ \h^*$ ,
$C \, = \, \sum_{i,j} C_{i,j} dx_i \o dx_j$ a closed meromorphic 2-form
on $\h^*$, \, $\epe$ a nonzero complex number. 
Then the function
$r \, : \, \h^* \, \to \g \oh  \g$ defined by
$$
r(\la) \, = \,  \sum_{i,j=1}^N  C_{i,j}(\la) x_i \otimes x_j \,
+ \,{\epe \over 2}\, \Omega \,+
\sum _{\al \in \Delta} \sum_{i=1}^{N_\al} \,
{\epe \over 2} \, \text{cotanh} \, ({\epe \over 2} \, (\al, \la -  \nu))\, 
e^i_\al \otimes  e^i_{-\al}\,
\tag 3.4
$$
is a classical dynamical r-matrix with the coupling constant $\epe$.
Here cotanh is the hyperbolic cotangent.}

\item {2.}{ Let $r \, : \, \h^* \, \to \g \oh \g$ be
a classical dynamical r-matrix with  nonzero coupling constant,
$r\,=\,r_\h\, +\, \sum_{\al \in \Dl}\,r_\al$ its weight
decomposition. 
If the function $r_\al$ is not constant for
any simple positive root $\al$, then the 
function $r$ has the form indicated in
(3.4).}
\endproclaim
Theorem 3.1 is proved in Section 3.6. 

For a simple Lie algebra $\g$,
the complete classification of  r-matrices 
$r \, : \, \h^* \, \to \g \oh \g$
with nonzero coupling constant
is given in Section 3.7.

\subhead 3.4. Classification of the classical dynamical r-matrices 
  with  zero coupling constant
\endsubhead

We shall use the notations of Theorem 3.1.

\proclaim {Theorem 3.2}

\item{1.}{ Let $X$ be a subset of the set of roots, $\Dl$, of a Kac-Moody
Lie algebra $\g$ such that  

1) if $\al,\beta\in X$ and $\al+\beta$ is a root then 
$\al+\beta\in X$, and

2) if $\al\in X$ then $-\al\in X$.
 
Let  $\nu$ be an element in $ \h^*$ ,
$C \, = \, \sum_{i,j} C_{i,j} dx_i \o dx_j$ a closed meromorphic 2-form
on $\h^*$. 
Then the function
$r \, : \, \h^* \, \to \g \oh \g$ defined by
$$
r(\la) \, = \,  \sum_{i,j=1}^N  C_{i,j}(\la) x_i \otimes x_j \,
+ \sum _{\al \in X} \sum_{i=1}^{N_\al} \,
{1\over (\al, \la -  \nu) } \,
e^i_\al \otimes  e^i_{-\al}\,
\tag 3.5
$$
is a classical dynamical r-matrix with zero coupling constant.
}
\item {2.}{ If $\g$ is a simple Lie algebra, then
any classical dynamical r-matrix with  zero
coupling constant has this form.}
\item {3.}{ Let $\g$ be an arbitrary Kac-Moody Lie algebra.
Let $r \, : \, \h^* \, \to \g \oh \g$ be
a classical dynamical r-matrix with  zero coupling constant,
$r\,=\,r_\h\, +\, \sum_{\al \in \Dl}\,r_\al$ its weight
decomposition. 
If the function $r_\al$ is not identically equal to zero for
any simple positive root $\al$, then the function $r$ has 
the form indicated in (3.5) with $X\,=\, \Dl$.}
\endproclaim

Theorem 3.2 is proved in Section 3.5. 

\subhead Example
\endsubhead
Let $\g$ be the Lie algebra of type $B_2$ with roots
$(\pm 1, 0), (0, \pm 1), (\pm 1, \pm 1)$. Then the set of long
roots $(\pm 1, \pm 1)$ gives an example of the set $X$. 

Consider an r-matrix of type (3.5). Assume that the element $\nu \in \h^*$
tends to infinity so that all terms of the matrix have limit.
Then the limiting function is an r-matrix of type (3.5) for a new set $X$.

Notice that the
r-matrix (3.5) corresponding to the example
is not a limiting case of the r-matrix
(3.5) with $X= \Dl$.

\subhead 3.5. Proof of Theorem 3.2 
\endsubhead

First we 
prove Theorem 3.2 assuming that $\g$ is a simple Lie algebra.
In this case dim $\g_\al = 1$ for any root $\al$. For any positive root
$\al$ fix elements $e_{\al} \in \g_{\al}$ and $e_{-\al} \in \g_{-\al}$
dual with respect to the bilinear form.

Let $r : \h^* \to \g \otimes \g$ be a classical dynamical r-matrix
with zero coupling constant.
The zero weight condition and the unitarity condition
imply that the r-matrix could be written
in the form
$$
r(\la) \, = \,\sum_{i,j=1}^N  C_{i,j}(\la) x_i \otimes x_j \,
+ \, \sum _{\al \in \Delta} \, \phi_\al (\la)\,
\,e_\al \otimes  e_{-\al}
$$
%
%
where $\phi_\al, \, C_{i,j}
 $ are suitable scalar meromorphic functions
such that
$\phi_{-\al}(\la) = - \phi_{\al}(\la)$, and
$C_{i,j}(\la) = - C_{j,i}(\la).$

The CDYB equation is an equation in $\g^{\otimes 3}$. 
The unitarity condition implies that the
left hand side of the CDYB equation is skew-symmetric with
respect to permutations of factors. 
This remark and the zero weight condition
show that in order to solve  the CDYB equation
it is enough to solve its
$\h \otimes \h \otimes \h-,$
 $\h \otimes \g_{\al} \otimes \g_{-\al}-$
 and $\g_\al \otimes \g_{\bt} \otimes \g_{\gm}-$parts,
where $ \al, \bt, \gm \in \Dl$ and in the last case
$\al +\bt + \gm =0$. 

A basis in 
$\g^{\otimes 3}$ is formed by the elements $x\otimes y\otimes z $
where $x, y, z$ run through $x_i, e_\al$, \, $i=1,...,N$,
$\al \in \Dl$.

The $x_i \otimes x_j \otimes x_k$-part  of the CDYB equation
has the form
$$
{\partial C_{j,k}\over \partial x_i}\,+\,
{\partial C_{k,i}\over \partial x_j}\,+\,
{\partial C_{i,j}\over \partial x_k}\,=\, 0
\tag 3.6
$$
and says that $\sum_{i,j=1}^N  C_{i,j}(\la) dx_i \o dx_j $
is a closed differential form.

The $\h\otimes \g_\al \otimes \g_{-\al}$-part of the CDYB equation
has the form
$$
\sum_k {\partial \phi_\al \over \partial x_k}\, x_k\otimes 
e_\al \otimes e_{-\al}\, + \, \phi_\al^2\, 
h_\al \otimes e_\al \otimes e_{-\al}\,=\,0
$$
where $h_\al=[e_\al, e_{-\al}]$.
This equation could be written in the form
$$
d\,\phi_\al\,+\,\phi_\al^2\,d\,h_\al =\,0.
\tag 3.7
$$
Hence $\phi_\al = 0$ or
 $\phi_\al = (h_\al - \nu_\al)^{-1}$ for some $\nu_\al \in \C$.
Here $h_\al $ is considered as a linear function on $\h^*$,
$h_\al(\la) = (\al, \la)$ for $\la \in \h^*$.

The $\g_\al\otimes \g_\bt \otimes \g_{-\al - \bt}$-part of the CDYB
equation has the form
$$
\phi_{\al}\, \phi_{\bt}\,-\, \phi_{\al}\, \phi_{\al +\bt}\,-\,
\phi_{\al +\bt}\, \phi_{\bt}\,=\,0
\tag 3.8
$$ 

\proclaim {Lemma 3.3}
Let $X=\{ \al \in \Dl \,\vert \, \phi_\al \neq 0 \}$. Then $X$ is closed
with respect to multiplication by $-1$ and addition.
\endproclaim
\demo {Proof}
The set $X$ is closed with respect to multiplication by $-1$ because of 
the unitarity condition. If $\phi_\al$ and $\phi_\bt$ are different from zero,
then $\phi_{\al +\bt}\, (\phi_\al + \phi_\bt )\,= \,\phi_\al \,\phi_\bt$,
and, hence, $\phi_{\al +\bt}$ is different from zero. $\square$
\enddemo
\proclaim {Lemma 3.4}
$$
\,\nu_{\al+\bt}\,=\,\nu_{\al}\,+\,\nu_{\bt}\, .
\tag 3.9
$$
\endproclaim
\demo{Proof}
Equation (3.8) implies
$$
{1\over h_{\al } - \nu_{\al} }\,{1\over h_{\bt} - \nu_{\bt} }\,=
( {1\over h_{\al} - \nu_{\al } }\,+\, {1\over h_{\bt} - \nu_{\bt} }\,)\,
{1\over h_{\al +\bt} - \nu_{\al +\bt} }\,.
$$
Since $h_{\al+\bt}=h_\al +h_\bt$,
in order to cancel the last pole one needs (3.9).
$\square$
\enddemo
\proclaim{Corollary 3.5} There is $\nu \in \h^*$ such that
$\nu_\al = (\al, \nu)$ for all $\al \in X$.
\endproclaim
This finishes  the proof of Theorem 3.2 for a simple Lie algebra.

Now we
assume that $\g$ is an arbitrary Kac-Moody algebra.
Let $r : \h^* \to \g \oh \g$ be a classical dynamical r-matrix
with  zero coupling constant.
The zero weight condition and the unitarity condition
imply that the r-matrix could be written
in the form
$$
r(\la) \, = \,  \sum_{i,j=1}^N  C_{i,j}(\la) x_i \otimes x_j \,
+ \sum _{\al \in \Dl} \sum_{i,j=1}^{N_\al} \,
\phi_\al^{i,j}(\la)\,
e^i_\al \otimes  e^j_{-\al}\,
$$
where $C_{i,j}(\la),\,\phi_\al^{i,j}(\la)\,$ are suitable scalar functions
such that $C_{i,j}=-C_{j,i}$ and $\phi_{-\al}^{j,i}=-
\phi_\al^{i,j}$.

The
$\h \otimes \h \otimes \h$-part  of the CDYB equation
is the same as for a simple Lie algebra and means that
 $\sum_{i,j=1}^N  C_{i,j}(\la) dx_i \o dx_j $
is a closed differential form.

To analyze the $\h \otimes \g_\al \otimes \g_{-\al}$-
and $\g_\al \otimes \g_{\bt} \otimes \g_{-\al - \bt}$-parts of 
the CDYB equation we shall introduce some useful linear operators
$\phi_\al(\la) : \g_\al \to \g_\al,$ where $\al \in \Dl$ and
$\la \in \h^*$. 
Namely, for any $\al \in \Dl$ and $\la \in \h^*$,
we define  $\phi_\al(\la)$  by the formula
$$
\phi_\al(\la) \,:\, e^j_\al \, \mapsto \sum_{i=1}^{N_\al}\,
\phi_\al^{i,j}(\la)\,e^i_\al \,.
\tag 3.10
$$
Let $r_\al(\la)$ denote the 
 $\g_\al \otimes \g_{-\al}$-part of the r-matrix,
$$
r_\al(\la)\,=\,
\sum_{i,j=1}^{N_\al}\,
\phi_\al^{i,j}(\la)\,e^i_\al \,\otimes e^j_{-\al}\,=
\,  \sum_{j=1}^{N_\al}\,
\phi_\al(\la)\,e^j_\al \,\otimes e^j_{-\al}\,.
$$

\proclaim{Lemma 3.6}
The $\h \otimes \g_\al \otimes \g_{-\al}$-part of the CDYB equation
has the form
$$
d\,\phi_\al\,+\,\phi_\al^2\,d\,h_\al =\,0.
\tag 3.11
$$
where the operator valued function
$\phi_\al : \h^* \to {\text{End}}\,(\g_\al)$
is defined in (3.10).
\endproclaim
\demo{Proof}
The  $\h \otimes \g_\al \otimes \g_{-\al}$-part of the CDYB equation
has the form
$$
\sum_{i,j=1}^{N_\al}\,\sum_k
{\partial \phi_\al^{i,j} \over \partial x_k}
(\la)\,x_k \otimes e^i_\al \,\otimes e^j_{-\al}\,
+\, [r_{-\al}^{12}(\la),\,r_{\al}^{13}(\la)]\,=\, 0.
$$
To prove the Lemma it is enough to show that
$$
[r_{-\al}^{12}(\la),\,r_{\al}^{13}(\la)]\,=
\sum_{l=1}^{N_\al}\,
h_\al \otimes (\phi_\al)^2e^l_\al \otimes
e^{l}_{-\al}\,.
$$
We  shall use the formula
$$
[e^i_{\al}, e^j_{-\al}]\,=\,(e^i_{\al}, e^j_{-\al})\,h_\al\,=\,\dl_{i,j}\,
h_\al \, .
$$
Then
$$
\align
[r_{-\al}^{12}(\la),\,r_{\al}^{13}(\la)]\,=\,
[\sum_{i,k}\,\phi^{i,k}_{-\al}\,e^{i}_{-\al}\otimes
e^{k}_{\al}\otimes 1,\,
\sum_{j,l}\,\phi^{j,l}_{\al}\,e^{j}_{\al}\otimes 1 \otimes
e^{l}_{-\al}\,]\,=
\\
\sum_{i,k,j,l}\,\phi^{i,k}_{-\al} \phi^{j,l}_{\al}\,
[e^i_{-\al},e^{j}_{\al}]\otimes e^k_\al \otimes
e^{l}_{-\al}\,=\,
\sum_{i,k,l}\,\phi^{k,i}_{\al} \phi^{i,l}_{\al}\,
h_\al \otimes e^k_\al \otimes
e^{l}_{-\al}\,=\,
\\
\sum_{l}\,
h_\al \otimes (\phi_\al)^2e^l_\al \otimes
e^{l}_{-\al}\,.
\endalign
$$
The Lemma is proved.
$\square$
\enddemo
Introduce a linear map
$$
\La : \g \otimes \g \otimes \g  \to \C,
\qquad 
x \otimes y \otimes z  \mapsto (x,[y,z])\,.
$$
Recall that the invariance of the  bilinear form implies
$(x,[y,z])=([x,y],z)$.

\proclaim{Lemma 3.7}

Let $\al, \bt,\gm \in \Dl$ be any roots such that
$\al+ \bt +\gm =0$. Then the $\g_{-\al} 
\otimes \g_{-\bt} \otimes \g_{-\gm}$-part
of the CDYB equation is equivalent to the statement that the composition
map
$$
\La \,\circ \,(\phi_\al \otimes \phi_\bt \otimes 1\,+\,
\phi_\al \otimes 1 \otimes \phi_\gm \,+\,
1 \otimes \phi_\bt \otimes \phi_\gm \,)\,:\,
\g_{\al} \otimes \g_{\bt} \otimes \g_{\gm}\, \to \, \C
\tag 3.12
$$
is the zero map.
\endproclaim
\demo{Proof} The $\g_{-\al} 
\otimes \g_{-\bt} \otimes \g_{-\gm}$-part
has the form
$$
[r_{-\al}^{13}(\la),\,r_{-\bt}^{23}(\la)]\,+\,
[r_{\bt}^{12}(\la),\,r_{\gm}^{13}(\la)]\,+\,
[r_{-\al}^{12}(\la),\,r_{\gm}^{23}(\la)]\,=\,0\,.
$$
Compute each of the terms.
$$
\align
&[r_{-\al}^{13}(\la),\,r_{-\bt}^{23}(\la)]\,=\,
[\sum_{i,j}\,\phi^{i,j}_{-\al}\,e^{i}_{-\al}\otimes 1 \otimes
e^{j}_{\al},\, 
\sum_{k,l}\,\phi^{k,l}_{-\bt}\,1 \otimes
e^{k}_{-\bt}\otimes e^{l}_{\bt}\,]\,=
\\
&\sum_{i,j,k,l}\,\phi^{i,j}_{-\al} \phi^{k,l}_{-\bt}\,
e^i_{-\al}\otimes e^{k}_{-\bt}\otimes [e^j_\al, e^l_\bt]\,=\,
\sum_{i,j,k,l,m}\,\phi^{i,j}_{-\al} \phi^{k,l}_{-\bt}\,
e^i_{-\al}\otimes e^{k}_{-\bt}\otimes e^m_{-\gm} \cdot
(e^m_{\gm}, [e^j_\al, e^l_\bt])\, .
\endalign
$$
Using the equalities $\sum_j \phi^{i,j}_{-\al} e^j_\al =
- \phi_\al  e^i_\al$,
$\sum_l \phi^{k,l}_{-\bt} e^l_\bt =
- \phi_\bt  e^k_\bt$, $(e^m_{\gm}, [e^j_\al, e^l_\bt]) = 
(e^j_{\al}, [e^l_\bt, e^m_\gm])$,
we get
$$
[r_{-\al}^{13}(\la),\,r_{-\bt}^{23}(\la)]\,=\,
\sum_{i,k,m} \, \La(\phi_\al e^i_\al \otimes \phi_\bt e^k_\bt
\otimes e^m_\gm ) \,
e^i_{-\al} \otimes e^k_{-\bt}
\otimes e^m_{-\gm} \,.
$$
Similarly
$$
\align
&[r_{\bt}^{12}(\la),\,r_{\gm}^{13}(\la)]\,=\,
\sum_{i,k,m} \, \La( e^i_\al \otimes \phi_\bt e^k_\bt
\otimes \phi_\gm e^m_\gm ) \, 
e^i_{-\al} \otimes e^k_{-\bt}
\otimes e^m_{-\gm} \,,
\\
&[r_{-\al}^{12}(\la),\,r_{\gm}^{13}(\la)]\,=\,
\sum_{i,k,m} \, \La(\phi_\al e^i_\al \otimes  e^k_\bt
\otimes \phi_\gm e^m_\gm ) \,
e^i_{-\al} \otimes e^k_{-\bt}
\otimes e^m_{-\gm} \,.
\endalign
$$
These formulae prove the Lemma.
$\square$
\enddemo
Lemmas 3.6 and 3.7 imply the first statement of Theorem 3.2 for an arbitrary 
Kac-Moody Lie algebra.

We shall show the third statement of Theorem 3.2 by induction on the
complexity of the roots. Let $\al_1, ..., \al_r$ be the simple positive roots
of $\g$ and $\al \in \Dl$ any
root. Represent it as a linear combination
of the simple roots, $\al=k_1 \al_1+...+k_r \al_r$. The absolute value
of the number $k_1 +...+k_r $ will be called the complexity of the root $\al$.
Thus the only roots of complexity one are the simple roots.

\proclaim{Lemma 3.8}

 Let $\g$ be a Kac-Moody Lie algebra. Let 
$r\,=\,r_\h\, +\, \sum_{\al \in \Dl}\,r_\al$
be a classical dynamical r-matrix with  zero
coupling constant and with nonzero components
$r_\al$ corresponding to the simple positive roots.
 Then the associated 
operator valued functions
$\phi_\al : \h^* \to \text{End}\, ( \g_\al)$
have the following form. There exists $\nu \in \h^*$ such that
$$
\phi_\al (\la) = {1\over (\al, \la-\nu)}\, \text{id} 
\tag 3.13
$$ 
for all $\al \in \Dl$ and $\la \in \h^*$.
\endproclaim
\demo{Proof } For any simple root $\al$, dim $\g_\al =1$.
Then $\phi_\al$ is a scalar function of the form
$\phi_\al(\la)=(h_\al - \nu_\al)^{-1}$ for some
$\nu_\al \in \C$, see Lemma 3.6. Hence there exists
 $\nu \in \h^*$ such that the operator
$\phi_\al(\la)$ has the form indicated in (3.13)
for all simple roots $\al$.
Our goal is to extend this formula to all roots.

Let $e_1,..., e_r , f_1,..., f_r$ be the Chevalley generators corresponding
to the simple positive
roots  $\al_1, ..., \al_r$.  For any root $\al$,
the space $\g_\al$ is generated by commutators of the elements
$e_1,..., e_r$, if $\al$ is positive, and by commutators of the elements
$ f_1,..., f_r$, if $\al$ is negative. Hence for any root $\al$
the space $\g_\al$ is generated by commutators of the form
$[x, y]$ where $x \in \g_\bt, \, y \in \g_\gm, \,$
$\al = \bt + \gm$, and the complexity of $\bt$ and $\gm$ is less than the
complexity of $\al$.

Assume now
that a root $\gm$ has the following property.
For any roots $\al$ and $ \bt$ such that $\al +\bt +\gm =0$ 
and the complexity of $\al$ and $\bt$ is less than the complexity
of $\gm$, the operators $\phi_\al$ and $\phi_\bt$ have the form indicated
in (3.13). Let us show that the operator $\phi_\gm$ has the same form.

Formula (3.12) implies
$$
({1 \over (\al, \la-\nu)}\,+\,{1\over (\bt, \la-\nu)})\,
([x,y],\, \phi_\gm(\la)\,z)\, = \, - \,
{1 \over (\al, \la-\nu)\, (\bt, \la-\nu)} ([x,y],\,z).
$$
This formula means that the operator 
$\phi_\gm(\la)$ is uniquely determined by the
operators $\phi_\al$ and $\phi_\bt$ corresponding to the roots
$\al$ and $\bt $ with the complexity less
than the complexity of $\gm$.
At the same time we see that if
the operator $\phi_\gm$ has the form indicated in (3.13),
then it satisfies both equations (3.12) and (3.11).
Lemma 3.8  and Theorem 3.2 are proved.
$\square$

\enddemo

\subhead 3.6. Proof of Theorem 3.1
\endsubhead

First we prove Theorem 3.1 assuming that $\g$ is a simple Lie algebra.
As in the proof of Theorem 3.2, fix the dual generators
$e_\al \in \g_\al$ and
$e_{-\al} \in \g_{-\al}$ for all roots $\al \in \Dl$.
Fix an orthonormal basis in $\h$, $x_1, ... , x_N$.

Let $r : \h^* \to \g \otimes \g$ be a meromorphic map, 
$\epe$ a nonzero complex number,
 $\Omega \in \g \otimes \g$ the Casimir operator of $\g$
associated to the bilinear form,
$\Om = \sum_k x_k\otimes x_k + \sum_{\al \in \Dl} e_\al
\otimes e_{-\al}$.

 Introduce a meromorphic map
$s : \h^* \to \g \otimes \g$
by the formula  $s(\la)\, = \, r(\la)\, - \, \epe \, \Omega / 2$.

\proclaim {Lemma 3.9} The map  
$r$ is  a classical dynamical r-matrix
with the coupling constant $\epe$ if and only
if the map $s$  satisfies the zero weight condition (3.1),
the unitarity condition,
$$
s^{12}(\la) + s^{21}(\la) = 0\, ,
$$
and the following analog of the CDYB equation,
$$
\align
\Alt (d s) \, + \, [s^{12},s^{13}]\,+ \, [s^{12},s^{23}]\,
+\,[s^{13},s^{23}] \, 
\tag 3.14
\\
+ \, {\epe^2 \over 4}\,(\, [\Omega^{12},
\Omega^{13}]\,+ \, [\Omega^{12},\Omega^{23}]\,
+\,[\Omega^{13},\Omega^{23}] \,)\, 
= \, 0 \, .
\endalign
$$
\endproclaim
\demo{Proof}
The only thing that needs to be checked is the fact 
that the terms of the form
$\epe \, [\Omega^{(i,j)}, s^{(k,l)}]/2$ 
cancel in the CDYB equation. This can be verified by
an easy direct calculation.
$\square$

\enddemo

The zero weight condition and the unitarity condition
imply that the matrix $s$ can be written
in the form
$$
s(\la) \, = \,\sum_{i,j=1}^N  C_{i,j}(\la) x_i \otimes x_j \,
+ \, \sum _{\al \in \Delta} \, \phi_\al (\la)\,
\,e_\al \otimes  e_{-\al}
$$
where $\phi_\al, \, C_{i,j}$ are suitable scalar meromorphic functions
such that
$\phi_{-\al}(\la) = - \phi_{\al}(\la)$, and
$C_{i,j}(\la) = - C_{j,i}(\la).$

The CDYB equation (3.14) is an equation in $\g^{\otimes 3}$.
Its left hand side is invariant with respect to even permutations of factors.
To solve the CDYB equation it suffices to solve its
$\h \otimes \h \otimes \h-,$
 $\h \otimes \g_{\al} \otimes \g_{-\al}-$
 and $\g_\al \otimes \g_{\bt} \otimes \g_{\gm}-$parts,
where $ \al, \bt, \gm \in \Dl$ and in the last case
$\al +\bt + \gm =0$.

The $x_i \otimes x_j \otimes x_k$-part  of the CDYB equation
has the form indicated in (3.6) and says that 
$\sum_{i,j=1}^N  C_{i,j}(\la) dx_i \o dx_j $
is a closed differential form.

The $\h\otimes \g_\al \otimes \g_{-\al}$-part of the CDYB equation
has the form
$$
\align
\sum_k {\partial \phi_\al \over \partial x_k}\, x_k\otimes 
e_\al \otimes e_{-\al}\, + \, \phi_\al^2\, 
h_\al \otimes e_\al \otimes e_{-\al}\,+\,
 {\epe^2\over 4}\, (\,-\,
h_\al \otimes e_\al \otimes e_{-\al}\,+&
\\
\sum_k \, (\, x_k\otimes [x_k, e_\al] \otimes e_{-\al}\, + \, 
x_k\otimes e_\al \otimes [x_k, e_{-\al}]&\,)) \, 
=\,0
\endalign
$$
where $h_\al=[e_\al, e_{-\al}]$.
This equation can be written in the form
$$
d\,\phi_\al\,+\,(\phi_\al^2\,-\,{\epe^2\over 4})\,d\,h_\al =\,0.
$$
Hence 
$$
\phi_\al(\la) \, = \,{\epe \over 2} \, 
\text{cotanh} \, ({\epe \over 2} \, ( h_\al - \nu_\al))\,  
$$
for some $\nu_\al \in \C$, or $\phi_\al\,=\,\pm \, \epe /2$.
Here $h_\al $ is considered as a linear function on $\h^*$,
cf. (3.7).

Let $\al, \bt, \gm \in \Dl$ be roots such that
$\al +\bt +\gm =0$.
The $\g_\al\otimes \g_\bt \otimes \g_{\gm}$-part of the CDYB
equation (3.14) has the form
$$
\phi_{\al}\, \phi_{\bt}\,+\, \phi_{\al}\, \phi_{\gm}\,+\,
\phi_{\gm}\, \phi_{\bt}
\,+\,{\epe^2\over 4}\,=\,0\,.
\tag  3.15
$$
If 
$$
\phi_\al(\la) \, = \,{\epe \over 2} \, 
\coth \, ({\epe \over 2} \, ( h_\al - \nu_\al))\,  ,
\qquad
\phi_\bt(\la) \, = \,{\epe \over 2} \, 
\coth \, ({\epe \over 2} \, ( h_\bt - \nu_\bt))\,  
$$
for some $\nu_\al, \nu_\bt \in \C$, then $\phi_\gm$ is not a constant,
so 
$$
\phi_\gm (\la) \, = \,{\epe \over 2} \, 
\coth \, ({\epe \over 2} \, ( h_\gm - \nu_\gm))\,  
$$
for some $\nu_\gm \in \C$. Starting from simple
positive roots we conclude that all functions $\phi_\al (\la)$
are not constant. Equation (3.15) also implies
$$
\coth \, ({\epe \over 2} \, (h_\gm - \nu_\gm))\,+
\, \coth \, ({\epe \over 2} \, (h_\al + h_\bt - \nu_\al -\nu_\bt))\,=\,0.
$$
Hence, there exists $\nu \in \h^*$ such that 
$\nu_\al = (\al, \nu)$ and
$$
\phi_\al(\la) \, = \,{\epe \over 2} \, 
\coth \, ({\epe \over 2} \, ( \al, \la - \nu))\,  
$$
for all roots $\al \in \Dl$.
Theorem 3.1 is proved for a simple Lie algebra $\g$.

The generalization of the above proof to the case of
a Kac-Moody algebra is word by word parallel
to the generalization of the proof of Theorem 3.2 from the simple Lie algebra
case to the Kac-Moody algebra case.

\subhead 3.7. Classification of classical dynamical r-matrices with nonzero coupling
constant, the simple Lie algebra case
\endsubhead

Let $\g $ be a simple Lie algebra, $\h$ its Cartan subalgebra,
$$
\g= \h \oplus \oplus _{\al \in \Delta} \g_\al 
$$
the root decomposition, $(\cdot,\cdot)$  an invariant nondegenerate
bilinear form on $\g$. For any positive root
$\al \in \Dl(\h)$ fix  basis elements $e_\al \in
\g_\al$ and $e_{-\al} \in \g_{-\al}$ which are dual with respect to
the bilinear form. Fix a basis  in the Cartan subalgebra,  $\{x_i\}$, 
orthonormal with respect to the bilinear form.

 Let 
$$
\Delta \,=\, \Delta_+ \, \cup \, \Delta_-
$$
be a polarization of roots into positive and negative,
$\Delta_+^s \subset \Delta_+ $ the set of simple positive roots.
Fix a subset $X \subset \Delta_+^s$ of the set of positive simple roots.

Fix a nonzero complex number $\epe$ and an element $\mu \,\in \, \h^*$.
For any root $\al$ introduce a meromorphic function $\phi_\al :
\h^* \to \C$ 
by the following rule. If a root $\al$ is a linear combination of
simple roots from $X$, then we set
$$
\phi_\al (\la) \,=\, 
{\epe \over 2} \, \coth \, ({\epe \over 2} \, (\al, \la -  \mu))\, .
$$
Otherwise we set $\phi_\al (\la) \,=\, \epe/2$, if $\al$ is positive,
and $\phi_\al (\la) \,=\, -\,\epe/2$, if $\al$ is negative.

Let $C \, = \, \sum_{i,j} C_{i,j} dx_i \o dx_j$ be
a closed meromorphic 2-form
on $\h^*$.

\proclaim {Theorem 3.10}

\item{1.}{ Introduce a function
 $r \, : \, \h^* \, \to \g \oh \g$ 
by the formula
$$
r(\la) \, = \,  \sum_{i,j=1}^N  C_{i,j}(\la) x_i \otimes x_j \,
+ \,{\epe \over 2}\, \Omega \,+
\sum _{\al \in \Delta} \,
\phi_\al (\la )\,
e_\al \otimes  e_{-\al}\,,
$$
where $C_{i,j}$ and $\phi_\al$ are defined above.
Then $r$ is a classical dynamical
r-matrix with nonzero coupling constant $\epe$. }
\item {2.}{ Any classical dynamical r-matrix, 
 $r \, : \, \h^* \, \to \g \otimes \g$ ,
with  nonzero coupling constant has this form. }
\endproclaim

The proof of Theorem 3.10 is based on the following  fact.

\proclaim {Theorem 3.11}

Let $Y \subset \Delta$ be a subset of the set of roots with two properties.

\item{A.}{ If $\al, \bt \in Y$ and $\al + \bt \in \Dl$, then
$\al + \bt \in Y$. }

\item{B.}{ If $\al$ is an element of $ Y$, then $-\al$ is not an element of $Y$. }

Then there exists a polarization
$
\Delta \,=\, \Delta_+ \, \cup \, \Delta_-
$
such that $Y \subset \Dl_+$.
\endproclaim
\demo{Proof of Theorem 3.11}
Consider 
$$
\frak n = \oplus_{\al \in Y}\, \C \, e_\al, \qquad
\frak m = \h\,\oplus \, \oplus_{\al \in Y}\, \C \, e_\al.
$$
Then $\frak n, \frak m $ are Lie subalgebras of $\g$.

\proclaim{Lemma 3.12}
\item{1.} { $
[\frak m , \frak m]\,=\, \frak n.
$
}
\item{2.} { 
The Killing form $B$ of $\frak m$ vanishes on $ \frak n$.
}
\endproclaim
\demo{Proof of Lemma 3.12} The first statement follows from (3.11.B).

The Killing form $B$ of $\frak m$ is defined by
$$
B(x,y)\,=\, tr\vert_{\frak m}\,(ad \,x\vert_{\frak m} \cdot ad \,y\vert_{\frak m})\,.
$$
Now the second statement follows from (3.11.B).
$\square$

\enddemo
According to Theorem 2.1.2 in [GG],
conditions 1 and 2 for a finite-dimensional Lie algebra $\frak m$
 imply that $\frak m$ is
solvable. Thus, $\frak m$ is a solvable subalgebra of $\g$.
This, in particular, means 
that $\frak m$ is contained in a Borel
subalgebra $\frak b$. Since $\h \subset \frak b$, the Borel subalgebra $\frak b$
defines a polarization of roots
$\Delta \,=\, \Delta_+ \, \cup \, \Delta_-$ such that $Y\subset \Dl_+$.
Theorem 3.11 is proved.
$\square$
\enddemo

\demo{Proof of  Theorem 3.10}
To prove the first statement of Theorem 3.10 it is enough to check that for
any roots $\al, \bt, \gm \in \Dl$ such that
$\al +\bt +\gm =0$, the functions $\phi_\al$, $\phi_\bt$,
$\phi_\gm$ satisfy equation (3.15). This could be easily done by direct
verification.

Now we prove the second statement of Theorem 3.10. Let 
 $r \, : \, \h^* \, \to \g \otimes \g$ be a classical dynamical r-matrix with
nonzero coupling constant $\epe$. According to Section 3.6, 
the r-matrix has the form
$$
r(\la) \, = \,  \sum_{i,j=1}^N  C_{i,j}(\la) x_i \otimes x_j \,
+ \,{\epe \over 2}\, \Omega \,+
\sum _{\al \in \Delta} \,
\phi_\al (\la )\,
e_\al \otimes  e_{-\al}\,,
$$
where  $\sum_{i,j} C_{i,j} dx_i \o dx_j$ is
a closed meromorphic 2-form on $\h^*$, and the functions
$\phi_\al$ are scalar meromorphic functions
such that $\phi_\al = {\epe \over 2}\, \coth \, ({\epe \over 2} ( (\al, \la) -  \mu_\al))$
for a suitable constant $\mu_\al$ or $\phi_\al = \pm \epe/2$. Moreover,
for any roots $\al, \bt, \gm \in \Dl$ such that
$\al +\bt +\gm =0$, the functions $\phi_\al$, $\phi_\bt$,
$\phi_\gm$ satisfy equation (3.15).

Let $Y \subset \Dl$ be the set of all roots $\al$ such that 
$\phi_\al = \epe/2$. Equation (3.15) and the unitarity condition 
easily imply that
the set $Y$ has properties 3.11.A - 3.11.B. 
By Theorem 3.11 there exists a polarization
of roots $\Delta \,=\, \Delta_+ \, \cup \, \Delta_-$ such that $Y\subset \Dl_+$.

Introduce two sets $X$ and $Z$ by
$$
 X\,=\, \Dl_+^s\, - \,\Dl_+^s \cap Y\,, \qquad
Z\,=\, \Dl_+ \,- Y.
$$
\proclaim{Lemma 3.13}

$Z$ is the span of $X$ in $\Dl_+$, i.e. $Z\,=\, \Z_{\geq 0}[X] \cap \Dl_+$.

\endproclaim
\demo{Proof}
If
$\al, \bt, \al + \bt \,\in \, \Dl_+$, $\phi_\al \neq \pm \epe/2$
and $\phi_\bt \neq \pm \epe/2$, then equation (3.15) implies that
$\phi_{\al+\bt} \neq \pm \epe/2$. This statement implies the inclusion
$ \Z_{\geq 0}[X] \,\subset \, Z$.

If $\al, \bt, \al + \bt \,\in \, \Dl_+$ and $\phi_\al = \epe/2$, 
then equation (3.15) implies that $\phi_{\al+\bt} = \epe/2$ ( since
$\phi_\bt$ could not be equal to $-\epe/2$ ).  This statement implies
the inclusion
$ Z\, \subset \, \Z_{\geq 0}[X]$. The Lemma is proved. $\square$

\enddemo

For $\al \in Z$ the functions $\phi_\al$ have the form
$\phi_\al = {\epe \over 2}\, \coth \, ({\epe \over 2} ( (\al, \la) -  \mu_\al))$
where $\mu_\al$ are suitable numbers. 
Moreover, if $\al, \bt, \al + \bt \,\in \, Z$ , then 
the corresponding constants $\mu_{\al}, \mu_{\bt}, \mu_{\al + \bt}$ 
satisfy the equation $\mu_{\al} + \mu_{\bt} = \mu_{\al + \bt}.$ 

Let $\mu \in \h^*$ be any element 
in the Cartan subalgebra such that
$\mu_\al \,=\, (\al, \mu)$ for all $\al \in X$. Then for all
 $\al \in Z $, we have $\mu_\al = (\al, \mu)$.

Now we may conclude that  the r-matrix $r$ has the form 
indicated in Theorem 3.10 and is associated to the polarization $\Dl = \Dl_+ 
\cup \Dl_-$, the set $X$ and the element $\mu$ constructed above.
Theorem 3.10 is proved. $\square$

\enddemo

Now we show that each of the r-matrices indicated in Theorem 3.10 is a limiting case
of the r-matrix (3.4). Namely, let $\Dl = \Dl_+ 
\cup \Dl_-$ be a polarization, $X \subset \Dl_+^s$ a subset of positive simple roots,
$\mu$ an element of $\h^*$, $\epe$ a nonzero complex number.
Set 
$$
\nu (t) \, =\, \mu\,+ \,t\, \sum_{i\in X}\, \omega_i
$$
where $\omega_i$ are fundamental weights of $\g$. Consider
the r-matrix $r_t(\la)$ defined by (3.4) with parameter $\nu$ equal to $\nu(t)$.
Then the limit of $r_t(\la)$ as $t$ tends to infinity is equal to the r-matrix
of Theorem 3.10 (unlike the example at the end of Section 3.4). 

\subhead Example
\endsubhead
For any polarization, the constant matrix
$$
r\,=\, {1\over 2}\,\sum_i \, x_i \otimes x_i
\,+\,\sum_{\al \in \Dl_+} \, e_\al \otimes e_{-\al}
\tag 3.16
$$
is a solution to the CDYB equation with coupling constant 1.
In particular, it is a solution to the classical Yang-Baxter equation.
This solution corresponds to $X=\Delta_+^s$. 

\subhead Remark
\endsubhead
Consider the r-matrix (3.4) with $C_{i,j}=0$ and $\nu=0$,
$$
r(\la) \, = \,  {\epe \over 2}\, \Omega \,+
\sum _{\al \in \Delta} \,
{\epe \over 2} \, \coth \, ({\epe \over 2} \, (\al, \la ))\, 
e_\al \otimes  e_{-\al}\,.
$$
Let $\Gamma \subset \h^*$ be an alcove of the Lie algebra $\g$
and $\Dl = \Dl_+ \cup \Dl_-$ the corresponding polarization.
If $\la$ tends to infinity inside the alcove $\Gamma$ in a generic direction,
then the r-matrix $r(\la)$ has a limit and this limit is given by
(3.16). Thus the r-matrix $r(\la)$ extrapolates different solutions 
of the classical Yang-Baxter equation of type
(3.16), labeled by different polarizations.

\subhead 3.8. Classical dynamical r-matrices associated with
a pair of Lie algebras
\endsubhead

Let $\g $ be a simple Lie algebra, $\h$ its Cartan subalgebra,
$$
\g= \h \oplus \oplus _{\al \in \Delta} \g_\al 
$$
the root decomposition, $(\cdot,\cdot)$  an invariant nondegenerate
bilinear form on $\g$ and $\Omega \in \g\otimes \g$
the associated Casimir operator.
Let $\l$ be a Lie subalgebra of $\g$ containing
$\h$. Assume that ${\frak l}$ is reductive. This condition is
equivalent to the condition  that 
there is a subset $\Dl(\l)_+ $ of the set
$\Dl_+$ of positive roots such that
$$
\l= \h \oplus \oplus _{\al \in \Delta(\l)_+} \, (\,
\g_\al \oplus \g_{-\al}\,) \,.
$$
Let $H \subset L \subset G$ be the corresponding complex Lie groups.

A meromorphic function
$$
r \, : \, \l^* \, \to \g \otimes \g
$$
is called {\it a 
classical dynamical r-matrix associated with the pair
$\l \subset \g$ }
if it satisfies the following
three conditions.
\item{1.}{{\it The invariance condition.} 
The meromorphic function $r$ is Ad $L$ - invariant,
$r(x\la x^{-1}) = $Ad $x \,
(\,r(\la)\,)$ for any $ \la \in \l^*$ and $x \in L$.}
\item{2.}{{\it The generalized unitarity,} 
$$
r^{12}(\la) + r^{21}(\la) = \epe \, \Omega
\tag 3.17
$$
for some constant $\epe \in \C$ and all $\la$.}
\item{3.}{{\it The classical dynamical Yang-Baxter equation, CDYBE,} 
$$
\Alt (d r) \, + \, [r^{12},r^{13}]\,+ \, [r^{12},r^{23}]\,
+\,[r^{13},r^{23}] \, = \, 0 \, .
\tag 3.18
$$
}

The differential of the r-matrix is considered in (3.18) as a 
meromorphic function
$$
dr\, : \, \l^* \, \to \g \otimes\g \otimes\g , \qquad
\lambda \, \mapsto \, \sum _i \, y_i \otimes {\partial r^{23} 
\over \partial y_i}(\la)\, , 
$$
where $\{y_i\}$ is any basis in $\l$. As before we 
denote by $\Alt (dr)$ the following
symmetrization of $dr$ ,
$$
\Alt (dr) \, = \, \sum _i \, y_i^{(1)} {\partial r^{23} 
\over \partial y_i}\, + \, \sum _i \, y_i^{(2)} {\partial r^{31} 
\over \partial y_i}\, + \, \sum _i \, y_i^{(3)} {\partial r^{12 }
\over \partial y_i}\, .
$$

It turns out that
the classification of the classical dynamical r-matrices associated 
with a pair $ \l \subset \g$  can be reduced to 
the classification of the  the classical dynamical r-matrices associated 
with the pair $\h \subset \g$ .

Namely, let us identify $\l^*$ with
$\h^* \oplus \oplus _{\al \in \Delta(\l)_+} \, (\,
\g_\al^* \oplus \g_{-\al}^*\,) \,.$
Let $r \, : \, \l^* \, \to \g \otimes \g$ be a 
 classical dynamical r-matrix associated 
with a pair $\l \subset \g$ .
First notice that if $\la \in \l^*$ is a semisimple element,
then there exists $x \in L$ such that $x\la x^{-1} \in \h^*$.
Since semisimple elements are dense  and since
the r-matrix satisfies 
the invariance condition, the function $r$ is completely determined
by its restriction to the dual to the
Cartan subalgebra, $r\vert_{\h^*}$ .

For any $\al \in \Dl(\l)_+$ fix  basis elements $e_\al \in
\g_\al$ and $e_{-\al} \in \g_{-\al}$ which are dual with respect to
the invariant bilinear form.
Define a function $\rho \, : \, \h^* \, \to \g \otimes \g$
by the formula
$$
\rho (\la)\,=\, \sum _{\al \in \Dl(\l)_+} \,
{e_\al \otimes  e_{-\al}\,- \, e_{-\al} \otimes  e_{\al}\,
\over (\al, \la ) } \, 
\tag 3.19
$$

\proclaim{ Theorem 3.14 }

A function $r \, : \, \l^* \, \to \g \otimes \g$ satisfying
the invariance condition is a classical dynamical r-matrix associated with
the pair $\l \subset \g$  if and only if
the function 
$$
r\vert_{\h^*} \,+\, \rho \,: \, \h^* \, \to \,\g \otimes \g
$$
is a classical dynamical r-matrix associated with the pair
$\h \subset \g$.

\endproclaim
Theorem 3.14 is proved in Section 3.9.

\subhead 3.9. Proof of Theorem 3.14
\endsubhead

Let  $\{x_i\}$ be any basis in $\h$. The elements
 $\{x_i\}$, \,and \, $e_\al, e_{-\al}, \, \al \in \Dl (\l)_+$,
form a basis in $\l$. Then
$$
dr(\la)\,= \, \sum _i \, x_i \otimes {\partial r^{23} 
\over \partial x_i}(\la)\, +
\sum _{\al \in \Dl (\l)_+}\, (\,
  e_\al \otimes {\partial r^{23} 
\over \partial e_\al}(\la)\, +\,
 e_{-\al} \otimes {\partial r^{23} 
\over \partial e_{-\al}}(\la)\,)\,
\tag 3.20
$$
for any $\la \in \l^*$.

\proclaim {Lemma 3.15}

Let $\la \in \h^*$. Then the second sum in (3.20) is equal to
$$
[\rho^{12}(\la)+\rho^{13}(\la),\, r^{23}(\la)].
$$
\endproclaim

\demo{Proof } 
 Let  $\{x_i^*\}$, \,and \, $e_\al^*, e_{-\al}^*, \, \al \in \Dl (\l)_+$
be the basis in $\l^*$ dual to the basis $\{x_i\}$, 
 \, $e_\al, e_{-\al}, \, \al \in \Dl (\l)_+$.
The Lie algebra $\l$ acts on $\l^*$ by ad$^*$.
Compute this action on elements of $\h^* \subset \l^*$.
For $\la \in \h^*$  we have
$((\text{ad}^*e_{\al})\,\la)(y) = - \la ([e_\al, y])$.
Hence $(\text{ad}^*e_{\al})\, \la = -\la (h_\al)\, e^*_{-\al} 
= -(\al, \la)\,e^*_{-\al} $. Similarly
$(\text{ad}^*e_{-\al})\, \la = 
 (\al , \la)\,e^*_{\al} $ and
$(\text{ad}^*x_j)\, \la = 0 $.

Now for any $\la \in \h^*$ we have
$$
{\partial r \over \partial e_{\al}}(\la)\,=\,
\,{1\over (\al, \la)}\,[e_{-\al} \otimes 1 \, +\,
1\otimes e_{-\al},\, r(\la)]\,.
\tag 3.21
$$
Indeed,
$$
r(e^{te_{-\al}} \la e^{-te_{-\al}})\,=\,r(\la) \,+\,
t\, (\al, \la)\, {\partial r \over \partial e_{\al}}(\la)
\,+\,O(t^2),
$$
and the invariance condition gives
$$
\align
r(e^{te_{-\al}} \la e^{-te_{-\al}})\,=\,
(e^{te_{-\al}}\otimes e^{te_{-\al}}) \, r(\la)\,
(e^{-te_{-\al}}\otimes e^{-te_{-\al}}) \,=\,
\\
r(\la)\,+\,t\,[e_{-\al} \otimes 1 \, +\,
1\otimes e_{-\al},\, r(\la)]\,+\,O(t^2)\,,
\endalign
$$
which implies (3.21).

Similarly
$$
{\partial r \over \partial e_{-\al}}(\la)\,=\,
\,{-\,1\over (\al, \la)}\,[e_{\al} \otimes 1 \, +\,
1\otimes e_{\al},\, r(\la)]\,.
\tag 3.22
$$
Formulae (3.21) and (3.22) prove the Lemma.
$\square$
\enddemo
Denote the first sum in (3.20) by $d_\h r$, then for $\la \in
\h^*$ we have
$$
dr(\la)\,=\,d_\h r(\la) \,+\,
[\rho^{12}(\la)+\rho^{13}(\la),\, r^{23}(\la)]\,.
$$

Now we finish the proof of Theorem 3.14. Let
$r \, : \, \l^* \, \to \g \otimes \g$ be a function
satisfying the invariance condition. Restrict this function to
$\h^*$. The invariance condition implies that the restriction satisfies
the zero weight condition (3.1). The restriction also satisfies
the generalized unitarity condition (3.17).
Introduce a function $\tilde{r} 
: \h^* \to \g \otimes \g$ by
$$
r\vert_{\h^*}(\la)\,=\,\tilde{r} (\la)\,-\,\rho (\la).
$$
Then the new function $\tilde{r} (\la)$ satisfies the zero weight condition
(3.1) and the generalized unitarity condition (3.17).

\proclaim{Lemma 3.16}

The function $r$ satisfies the CDYB equation (3.18)
on $\h^*$ if and only if the function $\tilde{r} $ satisfies the
CDYB equation (3.3).

\endproclaim
Lemma 3.16 implies Theorem 3.14.

\demo{Proof of Lemma 3.16}
The restriction of the CDYB equation on $\h^*$ takes the form
$$
\align
\Alt (d_\h \tilde{r} ) - \Alt (d_\h \rho) 
- [\rho^{12}+\rho^{13}, \rho^{23}] -
[\rho^{21}+\rho^{23}, \rho^{31}] 
\tag 3.23
\\
-[\rho^{31}+\rho^{32}, \rho^{12}] +
 [\rho^{12}, \rho^{13}]
+ [\rho^{12}, \rho^{23}]
+ [\rho^{13}, \rho^{23}] 
\\
-[\rho^{12}+\rho^{13}, \tilde{r}^{23}]
- [\rho^{21}+\rho^{23}, \tilde{r}^{31}]
- [\rho^{31}+\rho^{32}, \tilde{r}^{12}] 
\\
- [\tilde{r}^{12}, \rho^{13}] - [\rho^{12}, \tilde{r}^{13}]
- [\tilde{r}^{12}, \rho^{23}] - 
[\rho^{12}, \tilde{r}^{23}]
- [\tilde{r}^{13}, \rho^{23}]  
\\
-[\rho^{13}, \tilde{r}^{23}] +
 [\tilde{r}^{12}, \tilde{r}^{13}]
+ [\tilde{r}^{12}, \tilde{r}^{23}]
+ [\tilde{r}^{13}, \tilde{r}^{23}] \,.
\\
\endalign
$$
Using the identity 
$\rho^{12}(\la) + \rho^{21}(\la) = 0$
we conclude that the terms containing only the function
$\rho$ take the form of the CDYB equation (3.3).
By Theorem 3.2 we know that the function $\rho$ satisfies
equation (3.3), and hence the $\rho$-terms in (3.23)
are canceled out.

Using in addition the generalized unitarity condition
for $\tilde{r} $ we see that the $[\rho , \tilde{r} ]$-terms in (3.23)
are canceled out too. The remaining part of (3.23)
is the CDYB equation (3.3) for the function $\tilde{r} $.

Lemma 3.16 and Theorem 3.14 are proved.
$\square$
\enddemo

\head 4. Classification of classical dynamical $r$-matrices
with  spectral parameter
\endhead

\subhead  4.1. Classical dynamical r-matrices
with  spectral parameter
\endsubhead

Let $\g $ be a simple Lie algebra, $\h$ its Cartan subalgebra,
$$
\g= \h \oplus \oplus _{\al \in \Delta} \g_\al 
$$
the root decomposition, $(\cdot,\cdot)$  an invariant nondegenerate
bilinear form on $\g$  and $\Omega \in \g\otimes \g$
the associated Casimir operator.
 For any positive root
$\al \in \Dl(\h)$ fix  basis elements $e_\al \in
\g_\al$ and $e_{-\al} \in \g_{-\al}$ which are dual with respect to
the bilinear form. Fix a basis  in the Cartan subalgebra,  $\{x_i\}$, 
orthonormal with respect to the bilinear form.
Let $E$ be a neighborhood of $0$ in $\C$, $0 \in E \subset \C$.

A meromorphic function
$$
r \, : \, \h^* \,\times \,E \, \to \g \otimes \g
$$
is called {\it a classical
dynamical r-matrix with spectral parameter and associated with the 
pair $\h \subset \g$ }
if it satisfies the following
four conditions.
\item{1.}{{\it The zero weight condition,} 
$$
[h\otimes 1 + 1\otimes h\, , \, r(\lambda,z)]=0
\tag 4.1
$$
for all $\la \in \h^*, \, z \in E$ and $h \in \h$.}
\item{2.}{{\it The generalized unitarity,} 
$$
r^{12}(\la,z)\, + \,r^{21}(\la,-z) \,= \, 0\, .
\tag 4.2
$$
for  all $\la \in \h^*$ and $ z \in E$.}
\item{3.}{{\it The residue condition.}
$$
{\text{Res}}_{z=0} \, r(\la, z)\, = \, \epe \,\Omega\, . 
\tag 4.3
$$
for some constant $\epe \in \C$.
} 
\item{4.}{{\it The classical dynamical Yang-Baxter equation, CDYBE,} 
$$
\align
\Alt & (d_\h r) \, + \, [r^{12}(\la,z_{1,2}), r^{13}(\la,z_{1,3})]\,+ \, 
\tag 4.4
\\ &
[r^{12}(\la,z_{1,2}),r^{23}(\la,z_{2,3})]\,
+\,[r^{13}(\la,z_{1,3}),r^{23}(\la,z_{2,3})] \, = \, 0 \, .
\\
\endalign
$$
where $z_{i,j}=z_i-z_j$.
}

In (4.4) the differential of the r-matrix is considered  with respect to
the $\h$-variables, 
$$
d_\h r\, : \, \h^*\, \times E \, \to \g \otimes\g \otimes\g , \qquad
(\lambda , z) \, \mapsto \, \sum _i \, x_i \otimes {\partial r^{23} 
\over \partial x_i}(\la,z)\, .
$$
 In (4.4) we denote by $\Alt (d_\h r)$ the following
symmetrization of $d_\h r$ ,
$$
\Alt (d_\h r) \, = \, \sum _i \, x_i^{(1)} {\partial r^{23}
\over \partial x_i}(\la,z_{2,3})\, + \, \sum _i \, x_i^{(2)} {\partial r^{31} 
\over \partial x_i}(\la,z_{3,1})\, + \, \sum _i \, x_i^{(3)} {\partial r^{12 }
\over \partial x_i}(\la,z_{1,2})\, .
$$
The variable $z$ is called {\it the spectral parameter}.
The number $\epe$ in (4.3) is called {\it the coupling constant}.

We classify the germs of classical dynamical r-matrices with
spectral parameter at the subset $\h^* \times 0 \subset \h^* \times \C$.
We assume that an r-matrix has a Laurent power series expansion
of the form $r(\la,z) = \sum_m \, r^m (\la) \, z^m$ where
$r^m (\la)$ are meromorphic functions on $\h^*$. We assume that the Laurent
expansion is convergent to a meromorphic function of $\la$ and $z$ in a
punctured neighborhood of $\h^* \times 0 $ in $ \h^* \times \C$.
Any function $r(\la, z)$ with these properties will be called a function with
Laurent expansion.

\subhead
4.2. Gauge transformations
\endsubhead

In this subsection we introduce four  transformations of maps
$r \, : \, \h^* \times \C \, \to \g \otimes \g$ called {\it
the gauge transformations}. We assume that the map satisfies the
zero weight condition and the generalized unitarity condition
and, therefore,
has the form
$$
r(\la, z) \, = \,\sum_{i,j=1}^N  S_{i,j}(\la, z) x_i \otimes x_j \,
+ \, \sum _{\al \in \Delta} \, \phi_\al (\la, z)\,
\,e_\al \otimes  e_{-\al}
\tag 4.5
$$
where $\phi_\al, \, S_{i,j} $ are suitable scalar meromorphic functions
such that
$\phi_{-\al}(\la, -z) = - \phi_{\al}(\la, z)$, and
$S_{i,j}(\la, -z) = - S_{j,i}(\la, z).$

\item{1.}{  Let $C \, = \, \sum_{i,j} C_{i,j}(\la) dx_i \o dx_j$ be
a closed meromorphic 2-form. Set
$$
r(\la,z) \, \mapsto \, r(\la,z) \,+\,
\sum_{i,j=1}^N  C_{i,j}(\la)\, x_i \otimes x_j \,.
$$
}
\item{2.}{For a vector $v \in \h^*$  and a function $f$ on $\h^*$, we
denote  $L_vf$ a new function on $\h^*$ which is the derivative
of $f$ along the constant vector field defined by $v$.

For a holomorphic function 
$\psi \, : \, \h^* \, \to \, \C$, set 
$$
r(\la,z) \, \mapsto \,
\sum_{i,j=1}^N  \,(S_{i,j}(\la, z)\,+ \, z\, 
{\partial^2 \psi \over \partial x_i \, \partial x_j}
(\la))\,
x_i \otimes x_j \, +
\, \sum _{\al \in \Delta} \,\phi_\al (\la, z)\,
e^{z L_\al\psi (\la) }\,
e_{\al}\otimes e_{-\al} \,.
$$
}
\item{3.}{ For $ \nu \in \h^*$, set 
$$
r(\la,z) \, \mapsto \,  \sum_{i=1}^N  S_{i,j}(\la - \nu, z)\,
x_i \otimes x_j \, + \, \sum _{\al \in \Delta} \, \phi_\al (\la -\nu , z)\,
 \,e_\al \otimes  e_{-\al}\,.
$$
}
\item{4.}{ For nonzero complex number $a,\, b$, set 
$$
r(\la,z) \, \mapsto \,  a\,r(a \la , b z).
$$
}

Notice that the first three transformations do not change the residue of
$r(\la, z)$ at $z=0$ and the last transformation multiplies it by $a/b$.

\proclaim{Theorem 4.1}
Any gauge transformation transforms an r-matrix with spectral parameter
to an r-matrix with spectral parameter.
\endproclaim
Theorem 4.1 is proved in Section 4.5.

Two maps $r(\la, z)$ and  $r'(\la, z)$ will be called {\it equivalent}
if one of them could be transformed into another by a sequence
of gauge transformations.

\subhead
4.3. Elliptic, trigonometric and rational r-matrices with spectral parameter
\endsubhead

In this section we give examples of r-matrices with spectral parameter,
and  formulate a theorem that any r-matrix with
nonzero coupling constant is equivalent to one of the examples.

In order to describe  r-matrices we  use theta functions.
Let
$$
\theta_1(z,\tau) = - \sum^\infty_{j=-\infty} e^{\pi i(j+{1\over 2})^2\tau
+2\pi i (j+ {1\over 2})(z+{1\over 2})}
$$
be the Jacobi theta function. Let $\tau$ be a nonzero complex
number such that \newline
Im $\tau >0$. 
Following [FW] introduce the functions
$$
\sigma_w(z, \tau) = 
{\theta_1(w-z, \tau)\,\theta_1'(0, \tau)
\over
\theta_1(w, \tau)\,\theta_1(z, \tau)},
\qquad
\rho (z,  \tau) = 
{\theta_1'(z, \tau)
\over
\theta_1(z, \tau)}, 
\tag 4.6
$$
where $'$ means the derivative with respect to the first argument.
Notice that $\sigma_w(z, \tau) = z^{-1} + O(1)$, 
$\rho (z,  \tau) = z^{-1} + O(z)$ as $z \to 0$
and $\sigma_{-w}(- z, \tau) = - \sigma_w(z, \tau) $, 
$\rho ( - z,  \tau) = 
- \rho (z,  \tau)$.

\subhead  Example of an elliptic r-matrix
\endsubhead
$$
r(\la,z,\tau)\,=\, \rho(z,\tau)\sum_{i=1}^N x_i\otimes x_i\,
+ \,\sum_{\al \in \Dl} \sigma_{- (\al,\la)}(z,\tau) e_\al
\otimes e_{-\al}.
\tag 4.7
$$
 For every $\tau \in \C$, Im $\tau > 0$, the 
function $r(\la, z,\tau)$ is a classical dynamical r-matrix
with spectral parameter $z$ and  coupling constant $\epe= 1$ [FW].

\subhead  Examples of trigonometric r-matrices
\endsubhead

Let  $\Dl = \Dl_+ \cup \Dl_-$ be a polarization of the set of roots.
Fix a subset $X \subset \Dl_+^s$ of the set of simple positive roots.
For any root $\al$ introduce a meromorphic function $\phi_\al :
\h^* \to \C$ by the following rule. If a root $\al$ is a linear combination of
simple roots from $X$, then we set
$$
\phi_\al (\la, z) \,=\,{ \text{sin}\, ((\al, \la) + z)
\over
 \text{sin}\, (\al, \la) \,  \text{sin}\, z},
$$
otherwise we set
$$
\phi_\al (\la, z) \,=\,{ e^{-iz} 
\over
 \text{sin} \,z},
\qquad \text{for}\, \al \in \Dl_+, \qquad
\phi_\al (\la, z) \,=\,{ e^{iz} 
\over
 \text{sin}\, z},
\qquad \text{for} \, \al \in
 \Dl_-
  \,.
$$
We introduce a trigonometric r-matrix by
$$
r(\la,z)\,=\, \text{cotan}\, z \,\sum_{i=1}^N x_i\otimes x_i\,
+ \, \sum _{\al \in \Delta} \, \phi_\al (\la, z)\,
\,e_\al \otimes  e_{-\al}
\tag 4.8
$$
where $\text{cotan}\, z\,=\, \text{cos}\, z\,/ \text{sin} \,z.$

\subhead  Examples of rational r-matrices
\endsubhead

For a subset $X \subset \Dl$ of the set of roots
 closed with respect to addition and
multiplication by $-1$, we introduce a rational r-matrix by
$$
r(\la, z) \, = \,  {\Omega  \over z} \,
+ \sum _{\al \in X}  \,
{1\over (\al, \la) } \,
e_\al \otimes  e_{-\al}\, .
\tag 4.9
$$

\proclaim{Theorem 4.2} 
\item{1.} {Each of the r-matrices (4.7) - (4.9)
described in this section is a classical dynamical
 r-matrix with coupling constant $1$.}
\item{2.} {  The germ at
 $\h^* \times 0 \subset \h^* \times \C$ of any
 classical dynamical
r-matrix with spectral parameter
and a nonzero coupling constant is equivalent to one of the r-matrices
(4.7)- (4.9). }

\endproclaim

\proclaim{Corollary}  Any such a germ extends to a meromorphic function on
 $ \h^* \times \C$ .

\endproclaim

Theorem 4.2 is proved in Sections 4.4 and 4.5.

\subhead 4.4. Proof of part 1 of Theorem 4.2 
\endsubhead
According to 
[FW],
the elliptic r-matrix (4.7) is a classical dynamical
 r-matrix with spectral parameter .

Taking the limit of the r-matrix (4.7) when $\tau $ tends to $+ i \infty$,
we conclude that for any fixed element  $\nu \in \h$ the
r-matrix 
$$
r_\nu (\la,z)\,=\, \text{cotan}\, z \,\sum_{i=1}^N x_i\otimes x_i\,
+ \, \sum _{\al \in \Delta} \, 
\,{ \text{sin}\, ((\al, \la - \nu) + z)
\over
 \text{sin}\, (\al, \la - \nu) \,  \text{sin}\, z}\,
\,e_\al \otimes  e_{-\al}
$$
is a classical dynamical
 r-matrix with spectral parameter. If $\nu = 0$, then this r-matrix has
form  (4.8) with  $X = \Dl_+^s$. 
To show that any matrix of the form of (4.8) is a classical
dynamical r-matrix with spectral
parameter it is enough to apply to $r_\nu$
the limiting procedure with respect to $\nu$,
cf. the end of Section 3.7.

\proclaim{Lemma 4.3} If $r_0(\la)$ is a classical dynamical
 r-matrix without spectral parameter
and with zero coupling constant,
then 
$$
r(\la, z) \,=\, {\Omega \over z}\, + \, r_0(\la)
$$
is a classical dynamical
 r-matrix with spectral parameter.

\endproclaim
Proof by direct verification . $\square$

Lemma 4.3 and Theorem 3.2 show that any r-matrix (4.9)
is a classical dynamical r-matrix with spectral parameter.

\subhead 4.5. Proof of part 2 of Theorem 4.1  and Theorem 4.2 
\endsubhead

Let $r \, : \, \h^* \times \C \, \to \g \otimes \g$ be a germ at
 $\h^* \times 0 \subset \h^* \times \C$
 of a classical dynamical r-matrix with spectral parameter.
It follows from the zero weight condition that
the r-matrix can be written in the form of (4.5).

The CDYB equation (4.4)
is an equation in $\g^{\otimes 3}$. 
Its left hand side is invariant with respect to even permutations of factors
and simultaneuos permutations of variables $z_1, z_2, z_3$.
In order to  solve the CDYB equation it suffices to solve its
$\h \otimes \h \otimes \h-,$
 $\h \otimes \g_{\al} \otimes \g_{-\al}-$
 and $\g_\al \otimes \g_{\bt} \otimes \g_{\gm}-$parts,
where $ \al, \bt, \gm \in \Dl$ and in the last case
$\al +\bt + \gm =0$. 

{\bf 4.5.1. The $\h \otimes \h \otimes \h$-part of the CDYB equation.}

First we analyse the $\h \otimes \h \otimes \h$-part of the CDYB equation
(4.4) which has the form
$$
\sum _i \, x_i^{(1)} {\partial S^{23}
\over \partial x_i}(\la,z_{2}-z_{3}) +  \sum _i \, 
x_i^{(2)} {\partial S^{31} 
\over \partial x_i}(\la,z_{3}-z_1) +  \sum _i \, 
x_i^{(3)} {\partial S^{12 }
\over \partial x_i}(\la,z_{1}-z_2)=0.
\tag 4.10
$$

\proclaim{Lemma 4.4}

The sum
$\sum _i \, x_i^{(1)} {\partial S^{23}
\over \partial x_i}(\la,z)$
 is a linear function of $z$.

\endproclaim
\demo{Proof} Differentiating (4.10)
 with respect to $z_1$ and $z_2$ we conclude
that the second derivative $\partial^2 / 
\partial z_1 \partial z_2$ of the third
sum in (4.10) is equal to zero. This implies the Lemma.
$\square$

\enddemo
\proclaim {Corollary 4.5}

Consider the $\h \otimes \h$-valued function 
$S(\la, z) = \sum _{i,j}\, S_{i,j}(\la, z) x_i \otimes x_j$.
Let $S(\la, z)\,=\, \sum_n S^n(\la)\,z^n$ be its Laurent expansion.
Then
$$
{\partial S^n \over \partial x_i}(\la)\,=\,0
$$
for all $\la$, \,$i$, and $n, \,n \neq 0, 1$.

\endproclaim

Each of the three sums in (4.10)
is a linear function of $z_1, z_2, z_3$. Hence equation (4.10)
splits into
four independent equations corresponding to the coefficients of
 $z_1, z_2, z_3$ and the constant coefficient.
The constant coefficient part has the form
$$
{\partial S^0_{j,k}\over \partial x_i}\,+\,
{\partial S^0_{k,i}\over \partial x_j}\,+\,
{\partial S^0_{i,j}\over \partial x_k}\,=\, 0
$$
and is equivalent to the fact that
$\sum  S^0_{i,j}(\la) dx_i \o dx_j $
is a closed differential form.

\proclaim {Lemma 4.6}

There is a multivalued meromorphic function
 $\psi \, : \, \h^* \, \to \, \C$ 
with univalued meromorphic second derivatives
such that
$$
{\partial^2 \psi \over \partial x_i  \partial x_j}(\la)
\,=\, S^1_{i,j}(\la)
$$
for all $i, j, \la$.

\endproclaim
\demo{Proof }
The $z_1, z_2, z_3$-parts of equation (4.10)
together with the unitarity condition 
$S^1_{i,j}(\la)\,=\,S^1_{j,i}(\la)$
 have the form
$$
{\partial S^1_{i,j}\over \partial x_k}(\la)
\,=\,
{\partial S^1_{k,j}\over \partial x_i}(\la)
$$
for all $\la,
\,i, j, k$. These equations imply that there exist functions
$\phi_j$ such that $S^1_{i,j}= \partial \phi_j /\partial x_i$ and
moreover,  $\partial \phi_j /\partial x_i =
\partial \phi_i /\partial x_j$. Hence, there exists a function
$\psi $ with the properties indicated in the Lemma.
$\square$
\enddemo

{\bf Remark.} Later we will show that the function $\psi$ is in fact 
holomorphic in $\h^*$. 

\proclaim{Corollary 4.7}

Let $s\, : \, \C \, \to \, \h \otimes \h \,$ 
be a germ at $\h^* \times 0 \subset \h^* \times \C$ of a
meromorphic function
with Laurent expansion. Assume that 
$$
 s(z)\,+ \, s^{21}(-z) \,=\,0
$$
for all $z$, and
assume that the Laurent expansion of $s$
does not contain the terms of degree $0,\, 1$,
$s(z) = \sum_{m \neq  0, 1} \, s^m \, z^m$.
Let $r \, : \, \h^* \times \C \, \to \g \otimes \g$ be a germ at
$\h^* \times 0 \subset \h^* \times \C$ of a function
of the form
$$
\align
r(\la,z) \, = \, 
s(z)\,+\, \sum_{i,j=1}^N  C_{i,j}(\la)\, x_i \otimes x_j \,
+ \, z\, \sum_{i,j=1}^N  {\partial^2 \psi \over \partial x_i \, \partial x_j}
&(\la)\, x_i \otimes x_j \, +
\tag 4.11
\\
& \sum _{\al \in \Delta} \, \phi_\al (\la, z)\,
\,e_\al \otimes  e_{-\al} \,
\endalign
$$
where $\sum C_{i,j} dx_i \o dx_j$ is a closed meromorphic form on
$\h^*$, $\psi$ is a multivalued meromorphic function
 with univalued meromorphic second derivatives ,
 and the functions $\phi_\al$ are 
such that
$\phi_{-\al}(\la, -z) = - \phi_{\al}(\la, z)$. Then the function
$r$ satisfies the zero weight condition (4.1),
the unitarity condition (4.2)
 and the $\h \otimes \h \otimes \h$-part of the CDYB equation (4.4).
Moreover, any classical dynamical r-matrix with spectral
parameter has this form.

\endproclaim

{\bf 4.5.2. The $\h \otimes \g_{\al} \otimes \g_{-\al}$-part
of the CDYB equation.}

Now we analyze the $\h \otimes \g_{\al} \otimes \g_{-\al}$-part
of the CDYB equation. 
This part has the form
$$
\align
&\sum_{i}\, {\partial \phi_\al \over \partial x_i}\, 
x_i \otimes e_\al \otimes e_{-\al} \, +
\, \phi_{-\al}(\la, z_{1,2})\,\phi_{\al}(\la, z_{1,3})\,
[e_{-\al}, e_\al] \otimes e_\al \otimes e_{-\al} \, + \,
\\
&\sum_{i,j}\,\,\phi_{\al}(\la, z_{2,3})\,
(s_{i,j}(z_{1,2}) \,+\,C_{i,j}(\la)\,+\,z_{1,2}\,
 {\partial^2 \psi \over \partial x_i \, \partial x_j}
(\la))\,
x_i \otimes [x_j, e_\al] \otimes e_{-\al} \, +\,
\\
&\sum_{i,j}\,\,\phi_{\al}(\la, z_{2,3})\,
(s_{i,j}(z_{1,3}) \,+\,C_{i,j}(\la)\,+\,z_{1,3}\,
 {\partial^2 \psi \over \partial x_i \, \partial x_j}
(\la))\,
x_i \otimes e_\al \otimes [x_j, e_{-\al}] \, =\,0 .
\endalign
$$
This equation can be written as
$$
\align
-\phi_{-\al}(\la, z_{1,2})\,\phi_{\al}(\la, z_{1,3})\,
h_\al \otimes e_\al \otimes e_{-\al} \, +
\tag 4.12
\\
\sum_{i}\,\lbrace(\, {\partial \phi_\al \over \partial x_i}\, 
+ \, \\
\phi_{\al}(\la, z_{2,3})\,
\sum_j \, [\,
s_{i,j}(z_{1,2}) \,-\,
s_{i,j}(z_{1,3}) \,+\,z_{3,2}\,
\al(x_j)\,
 {\partial^2 \psi \over \partial x_i \, \partial x_j}
(\la)\,]\,\rbrace)\,
x_i \otimes e_\al \otimes e_{-\al} \, =\,0.
\endalign
$$
We interpret this equation as an equation of differential
1-forms on $\h^*$
identifying linear functions with their differentials
and ignoring the factor $ e_\al \otimes e_{-\al}$.
Then the first term in (4.12)
can be written as 
$ \phi_\al(\la, z_{2,1})\, \phi_\al(\la, z_{1,3})\,dh_\al$.
The second has the form $d \phi_\al(\la, z_{2,3})$
where the differential is with respect to $\la$. 
For a fixed $z$ consider $\sum s_{i,j}(z)\, x_i\otimes x_j$ as a
bilinear form on $\h^*$,
$$
s(z)\{\la, \mu\}=\sum s_{i,j}(z) <x_i, \la> <x_j, \mu>.
$$ 
If the second argument of this form is equal to $\al$, then we get a
linear function on $\h^*$,
$\sum  s_{i,j}(z) <x_j, \al> x_i.$  Hence the third and the forth
terms in (4.12) have the form
$ \phi_{\al}(\la, z_{2,3})\, (
ds(z_{1,2})\{\la, \al\} \,-\,ds(z_{1,3})\{\la, \al\})$
where  the differentials are with respect to $\la$.
Finally the last term in (4.12) could be written as
$z_{2,3}\phi_{\al}(\la, z_{2,3})
\, L_\al d \psi (\la)$ where $L_\al d \psi (\la)$ 
is the Lie
derivative of the differential with respect to the constant vector
field defined by $\al$. Now the 
$\h \otimes \g_{\al} \otimes \g_{-\al}$-part
of the CDYB equation takes the form
$$
\align
 d \phi_\al(\la, z_{2,3})\,&      +\,
\phi_\al(\la, z_{2,1})\, \phi_\al(\la, z_{1,3})\,dh_\al \,+
\\
& 
\phi_{\al}(\la, z_{2,3})\, (
ds(z_{1,2})\{\la, \al\} \,-\,d s(z_{1,3})\{\la, \al\}\,+
\,z_{2,3}\, L_\al d \psi (\la)\,)\,=\,0.
\\
\endalign
$$
Setting $u=z_{2,1}$ and $v=z_{1,3}$
we get
$$
\align
{d \phi_\al(\la, u+v) \over \phi_\al(\la, u+v)}
\,  &  
+\,
{ \phi_\al(\la, u)\, \phi_\al(\la, v) \over
 \phi_\al(\la, u+v)}\,dh_\al \,+\,
\\
ds &
(-u)\{\la, \al\} \,-\,d s(v)\{\la, \al\}\,-
\,(u+v)\, L_\al d \psi (\la)\,
\,=\,0
\endalign
$$
where the differentials are with respect to $\la$.

Make a change of variables, 
$\phi_\al (\la, z) = \Phi (\la, z)\, e^{z L_\al \psi(\la)}$.
Then
$$
{d \Phi(\la, u+v) \over \Phi(\la, u+v)}
\,+\,
{ \Phi (\la, u)\, \Phi(\la, v) \over
 \Phi(\la, u+v)}\,dh_\al \,+\,
ds(-u)\{\la, \al\} \,-\,d s(v)\{\la, \al\}\,
\,=\,0.
$$
Taking the 
differential of both sides we see that the second ratio depends 
only on $h_\al$. Hence the function $\Phi$ has the form
$$
\Phi (\la, z)\,=\, \mu (h_\al (\la), z)\, e^{z \nu (\la)}
$$
for suitable new functions 
$\mu (h_\al (\la), z)$ and $\nu (\la)$.
Now the equation takes the form
$$
\align
{ \partial \mu(h_\al, u+v) / \partial h_\al
 \over \mu(h_\al, u+v)}\, dh_\al
\,+\, (u+v)\, d \nu & (\la)\,+\,
{ \mu (h_\al, u)\, \mu (h_\al, v) \over
 \mu (h_\al, u+v)}\,dh_\al \,+\,
\tag 4.13
\\
& ds(-u)\{\la, \al\} \,-\,d s(v)\{\la, \al\}\,
\,=\,0.
\endalign
$$
Consider a coordinate system on $\h^*$, \
$y_1, ..., y_N \in \h$, such that $y_1 = h_\al$
and 
\newline
$<\al, y_i>=0$ for $i>1.$ 
Then (4.13) gives the equations
$$
\align
{ \partial \mu(h_\al, u+v) / \partial h_\al
 \over \mu(h_\al, u+v)}\, 
\,+\, (u+v)\,  {\partial \nu\over \partial h_\al}& (\la)\,+\,
{ \mu (h_\al, u)\, \mu (h_\al, v) \over
 \mu (h_\al, u+v)} \,+\,
\tag 4.14
\\
&
{s(-u)\{\al, \al\} \,-\,s(v)\{\al, \al\} \over (\al, \al)}
\,
\,=\,0,
\endalign
$$
$$
(u+v)\, {\partial \nu\over \partial y_i} (\la)\,+\,
s(-u)\{{\partial \over \partial y_i}
, \al\} \,-\, s(v)\{{\partial \over \partial y_i}
, \al\}\,
\,=\,0 \, , \qquad i=2,...,N\, .
\tag 4.15
$$
(Here $\{\frac{\partial}{\partial y_i}\}$ is the basis of $\h^*$ 
dual to the basis $\{y_i\}$ of $\h$.)
Equations (4.15) imply that
$\partial \nu /\partial y_i =0$ since the Laurent expansion of
the function $s(z)$ does not have the first order term.
So we can assume that $\nu = 0$ and, therefore,
$\phi_\al (\la, z)\, = \, \mu (h_\al (\la), z)\, e^{z L_\al \psi(\la)}$.

Now equations (4.15) imply that 
$s(z)\{{\partial \over \partial y_i},\al\}$ does not depend on
$z$. But the Laurent expansion of this function does not contain
the terms of degree zero and one. So
the function $s(z)\{{\partial \over \partial y_i},\al\}$ is identically
equal to zero for $i \geq 2$.

The fact that $s(z)\{{\partial \over \partial y_i},\al\}$ is identically
equal to zero for $i \geq 2$ easily implies that 
$$
s(z)\,= \, t(z)\,\sum_{i=1}^N\, x_i \otimes x_i\,
$$
where $t(z)$ 
is a scalar function. The function $t(z)$ is a nonzero function since
the coupling constant is not zero.
We have $t(-z)= - t(z)$, and the Laurent expansion of $t$ does not
have the terms of degree zero and one.
Then (4.14) takes the form
$$
{ \partial \mu(h_\al, u+v) \over \partial h_\al}\,=\,
 (t(u) \,+\,t(v)) \mu(h_\al, u+v)\, 
-\, 
\mu (h_\al,  u)\, \mu (h_\al, v) \,.
\tag 4.16
$$

{\bf 4.5.3. Proof of Theorem 4.1.} 

Before proceeding with analysis of equation (4.16) let us write
the $\g_\al \otimes \g_{\bt} \otimes \g_{\gm}-$ part of the CDYB equation,
$$
\phi_{\al}(\la, z_{13})\,\phi_{\bt}(\la, z_{23})\,+\,
\phi_{\bt}(\la, z_{21})\,\phi_{\gm}(\la, z_{31})\,+\,
\phi_{\al}(\la, z_{12})\,\phi_{\gm}(\la, z_{32})\,=\,0\,.
\tag 4.17
$$
It is easy to see that equations (4.10),(4.16), (4.17) are invariant 
with respect to the gauge transformations of Section 4.2. This proves 
Theorem 4.1.

{\bf 4.5.4. The $\h \otimes \g_{\al} \otimes \g_{-\al}$-part
of the CDYB equation: classification of solutions.}

Now we will find all solutions of equation (4.16).

\proclaim{Lemma 4.8}

The function $\mu (x, z)$ has at most a simple pole at $z=0$.

\endproclaim

\demo{Proof}
 Applying the operator $\partial /\partial u
- \partial /\partial v$ to both sides of (4.16) we have
$$
(t'(u) \,-\,t'(v)) \mu(x, u+v)\, 
-\, 
\mu '(x,  u)\, \mu (x, v) \,+\,
\mu (x,  u)\, \mu '(x, v) 
\,=\, 0.
$$
Set $v=-u+\dl$. Then using the fact that $t'(-z)=t'(z)$
we get
$$
(t'(u) \,-\,t'(u-\dl)) \mu(x, \dl)\, 
-\, 
\mu '(x,  u)\, \mu (x, -u+\dl) \,+\,
\mu (x,  u)\, \mu '(x, -u+\dl) 
\,=\, 0.
$$
The function $t(z)$ could not be a linear function. Hence 
the first factor is of order $\dl$.
The second and the third terms are regular at generic values of $u$.
This shows the Lemma.
$\square$
\enddemo

The Lemma easily implies that the function $t(z)$ has also
 at most a simple pole at $z=0$. Thus,
$$
\mu (x, z) = \sum_{n=-1}^\infty \mu_n(x) z^n, \qquad
 t(z) = \sum_{n=-1,\, n \,\text{odd}}^\infty t_n z^n .
$$
Notice that from the residue condition (4.3) we know that 
$
t_{-1} = \epe \, (\al, \al)
$
where $\epe$ is the coupling constant.
Substitute the expansions
 into (4.16) and multiply both sides by $uv(u+v)$,
$$
\align
\sum_{n=-1}^\infty \mu'_n(x) (u+v)^{n+1}uv -
\sum_{n,m=-1}^\infty t_n \mu_m(x)(u+v)^{m+1}(u^n+v^n)uv+&
\tag 4.18
\\
\sum_{n,m=-1}^\infty \mu_n(x)\mu_m(x)u^{n+1}v^{m+1}&
(u+v)
=0.
\endalign
$$
This is an equality of power series in $u$ and $v$.
Equating the homogeneous parts we get a sequence of
equations for the numbers $\{t_n\}$ and the functions $\{\mu_n(x)\}$. 

We write the first equations. The equation of degree $1$ has the form
$$
\mu_{-1}(x)\,=\,t_{-1}\,.
\tag 4.19
$$
The equation of degree $2$ has the form $\mu_{-1}'(x)=0$
and follows from (4.19). 
The equations of degree $3$ and $4$ have the form 
$$
\align
&\mu_{0}'\,=\,2t_{-1}\mu_1\,-\,\mu_0^2\,+\,t_{-1}t_{1}\,,
\\
&\mu_{1}'\,=\,3t_{-1}\mu_{2}\,+\,t_{1}\mu_0\,-\,\mu_0\mu_1\,.
\endalign
$$
The equation of degree $5$ gives two scalar equations 
$$
\align
&\mu_{2}'\,=\,4t_{-1}\mu_{3}\,+\,t_{1}\mu_1\,+
\,t_3t_{-1}-\mu_0\mu_2\,,
\\
&2\mu_{2}'\,=\,6t_{-1}\mu_{3}\,+\,2t_{1}\mu_1\,-
\,t_3t_{-1}-\mu_1^2\,,
\endalign
$$
which imply
$$
\align
&\mu_{2}'\,=\,t_{1}\mu_{1}\, -\,5t_{-1}t_{3}\,+\,3\mu_0\mu_2\,
-\,2\mu_1^2,
\\
& 2\mu_3t_{-1}\,=\,2\mu_0\mu_2\,-\,3t_{-1}t_3\,-\,\mu_1^2.
\endalign
$$
Thus we have
$$
\align
&\mu_{0}'\,=\,t_{-1}t_{1}\,+\,2t_{-1}\mu_1\,-\,\mu_0^2\,,
\tag 4.20
\\
&\mu_{1}'\,=\,3t_{-1}\mu_{2}\,+\,t_{1}\mu_0\,-\,\mu_0\mu_1\,,
\\
&\mu_{2}'\,=\,t_{1}\mu_{1}\, -\,5t_{-1}t_{3}\,+\,3\mu_0\mu_2\,
-\,2\mu_1^2.
\endalign
$$

\proclaim{Lemma 4.9}

Let $t_{-1}\neq 0$ and $n\geq 3$. Then the degree $n+3$ equation
determine $\mu_{n+1}, t_{n+1}$ uniquely in terms of
$\mu_{m}, t_{m}$ with smaller $m$.

\endproclaim

\demo{Proof}
The degree $n+3$ equation contains only $\mu_{m}, t_{m}$
with $m\leq n+1$ and has the form
$$
\align
- t_{-1}\mu_{n+1}(u+v)^{n+2}({1\over u}+
{1\over v}) u v\, - \,
&t_{n+1}\mu_{-1}(u^{n+1}+v^{n+1}) u v\, +
\\
& \mu_{-1} \mu_{n+1}(u^{n+2}+v^{n+2}) (u + v) +\, ... \,=\,0
\\
\endalign
$$
where $...$ denotes the terms containing only $t_m$ and $\mu_m$
with $m<n+1$.

The coefficients of $u^{n+2}v$ and  $u^{n+1}v^2$ have the form
$$
\align
&(n+2) \mu_{n+1} + t_{n+1}  + ...  = 0,
\\
&(n+3)(n+2) \mu_{n+1} + ...  = 0.
\endalign
$$
The equations imply the Lemma.
$\square$ \enddemo
\proclaim{Corollary 4.10} If
 $t_{-1}\neq 0$, then a solution of equation (4.16)
is uniquely determined by the six parameters
$\mu_0(0), \mu_1(0), \mu_2(0), t_{-1}, t_{1}, t_3$.

\endproclaim

Now we present a six parameter family of solutions.
We shall use the solution of the CDYB equation (4.4)
given in [FW],
$$
r(\la,z,\tau)\,=\, -\rho(z,\tau)\sum_i x_i\otimes x_i\,
-\,\sum_{\al \in \Dl} \sigma_{h_\al(\la)}(z,\tau) e_\al
\otimes e_{-\al}
$$
where the functions $\rho$ and $\sigma$  are defined
in (4.6). For every $\tau \in \C$, Im $\tau > 0$, the 
function $r(\la, z,\tau)$ is a classical dynamical r-matrix
with spectral parameter $z$ and  coupling constant $\epe=-1$ [FW].
Hence, for every  $\tau$,
the functions $t(z)=-\rho(z,\tau)$
and $\mu(x, z)=-\sigma_{x}(z,\tau)$
form a solution of (4.16).
\proclaim{Lemma 4.11}
Let $t(z)$ and $\mu(x, z)$ be a solution of (4.16).
Let $A$ be a complex number, then
$$
\align
t(Az), &\qquad \mu(x, Az),
\\
t(z), &\qquad \mu(x + A, z),
\\
A\,t(z), &\qquad A\,\mu(Ax, z),
\\
t(z) + Az, &\qquad e^{A z x}\,\mu(x, z),
\\
t(z), &\qquad e^{A z }\,\mu(x, z)
\endalign
$$
are solutions of (4.16).
\endproclaim
\proclaim{Corollary 4.12}
The functions 
$$
\align
& t(z)=-A\rho(B z,\tau)+D z, 
\tag 4.21
\\
& \mu(x, z)=-A\sigma_{A x-C}
(B z,\tau)
e^{z(D x+E)}
\endalign
$$
form a solution of (4.16) depending on six parameters
$A, B, C, D, E, \tau$.

\endproclaim
Let $t, \mu$ be the solution of (4.16) defined by (4.21).
Let $t_{-1}, t_1, t_3,$
$ \mu_{0}(0),$ $ \mu_{1}(0), $ $
\mu_{2}(0)$ be the corresponding Taylor 
coefficients of $t$ and $\mu$, cf. Corollary 4.10.
These six Taylor coefficients are functions of the six parameters
$A, B, C, D, E, \tau$ and define a meromorphic 
map $\chi: \C^5 \times H \to \C^6$,
where $H$ is the upper half plane.

\proclaim{Lemma 4.13}
The Jacobian of this map 
is not identically equal to zero, and
the image of the map is dense in $\C^6$.

\endproclaim
\demo {Proof}
In order to prove that the Jacobian is not zero
it suffices to show that all the six parameters $A,...,\tau$ could be 
recovered from 
the function
$$
\mu=-A{\theta_1(Ax-C-Bz, \tau) \theta_1'(0, \tau)
\over \theta_1(Ax-C, \tau) \theta_1(Bz, \tau) }
e^{z(Dx + E)}.
$$
In fact, the poles of this function are given by
$$
Ax-C \in \Z \oplus \tau \Z \qquad
\text{and} \qquad Bz \in \Z \oplus \tau \Z .
$$
Knowing the poles we recover $\tau, A, B, C$.
If $Bz \mapsto Bz + 1$, then
$\mu \mapsto \mu e^{(Dx + E)/B}$,
this property allows us to recover $D, E$ at least locally.

Now let us show the density of the image of $\chi$. 
Let $p(z)=\frac{d}{d z}\ln \mu(0,z)=-\frac{1}{z}+p_0+p_1z+p_2z^2+...$. 
Consider the map $\tilde\chi: \Bbb C^5\times H\to \Bbb C^6$ 
given by the formula $(A,B,D,E,C,\tau)\to (t_{-1},t_1,t_3,p_0,p_1,p_2)$. 
It is enough to show that the image of $\tilde\chi$ is dense. 

We have
$$
t(z)=-A\rho(Bz,\tau)+Dz, \qquad 
p(z)=E+\frac{d}{dz}\ln\frac{\theta_1(C+Bz,\tau)}
{\theta_1(Bz,\tau)}.
$$
Thus, $t_1=D+f_1(A,B,\tau)$, $p_0=E+f_2(B, C, \tau)$, for some meromorphic 
functions $f_1,f_2$, and $t_{-1},t_3,p_1,p_2$ do not depend on $D$ and $E$. 
Thus, to show the density of the image of $\tilde\chi$, 
it is enough to show the density of the image of 
$\xi: \C^3\times H\to \C^4$ given by 
$(A,B,C,\tau)\to (t_{-1},t_3,p_1,p_2)$.

It is clear that the function $p'(z)$ 
is doubly periodic with respect to $C$, 
with periods $1,\tau$. Therefore, for any fixed $\tau$, the map 
$\xi$ is a {\it rational} map $\xi_\tau: 
\Bbb C^2\times E_\tau\to \Bbb C^4$, where $E_\tau$ is the elliptic curve 
corresponding to $\tau$. Denote by $I_\tau$ 
the closure of the image of $\xi_\tau$. 
As the Jacobian of $\xi$ is not identically zero, for generic $\tau$, the set
$I_\tau$ is an irreducible algebraic hypersurface in $\Bbb C^4$.

It is easy to see that each of the functions
 $t''(z),p'(z)$ satisfies the modular 
invariance conditions  
$$
\gather
f(A,B,C,z,\tau+2)=f(A,B,C,z,\tau),
\\
f(A,B,C,z,-1/\tau)=
f(A\tau,B\tau,C\tau,z,\tau).
\endgather
$$
Therefore, the hypersurface $I_\tau$ is modular invariant
with respect to the subgroup $\tilde\Gamma$ of the modular group $\Gamma$ 
generated
by $\tau\to \tau+2$, and $\tau\to -1/\tau$. 
This means that the coefficients of the equation of $I_\tau$ are 
modular functions 
on the modular curve $\Sigma=H/\tilde\Gamma$
(it is easy to see that they have power growth 
in $q$ as $q=e^{2\pi i\tau}\to 0$).
This shows that the hypersurfaces $I_\tau$ form an algebraic family
over $\Sigma$. Let $T$ 
be the total space of this family. We have a natural rational 
map $\psi: T\to \C^4$, and the closure $I$ of its image
coincides with the closure of the image of $\xi$. 
Since the map is rational and has a nonzero Jacobian, we have 
$I=\Bbb C^4$, as desired.  
$\square$
\enddemo

Corollaries 4.10, \, 4.12 and Lemma 4.13 tell us that  all solutions
$t, \, \mu$ to equation (4.16) are limits of solutions given by (4.21).
It is enough to list the limits of the function $\mu$ since then
the function $t$ can be recovered from (4.16).
Without loss of generality we assume that  the coupling constant
$\epe$ is equal to $1$, i.e. $A = -B$.

Let 
$$
f(x, z) \,=\, B\, 
{\theta (B (x - x_0 -z ), \tau) \, \theta ' (0, \tau) 
\over
\theta (B (x - x_0 ), \tau) \, \theta  ( B z , \tau) }
\, e^{z(D x + E)}
$$
be a function of $x, z $ depending on parameters $B, D, E, x_0, \tau $.
A function $g(x,z)$ will be called a limit of the function $f$ if there exist
sequences $B_n, D_n, E_n, x_{0,n}, \tau_n$ such that $g(x,z)$ is the limit of
$f(x, z ; B_n, D_n, E_n, x_{0,n}, \tau_n)$ when $n$ tends to infinity.

\proclaim{Proposition 4.14}

Any limit $g$ of the function $f$ has one of the following three forms.

Rational type.
$$
\align
g\,=& \,{1 \over z}\,e^{z(D x + E)}\,, \tag 4.22
\\
g\,=& \,{ x - x_0 - z \over (x - x_0) z }\,e^{z(D x + E)}\, \tag 4.23
\endalign
$$
where $ D, E, x_0$ are parameters.

Trigonometric type.
$$
\align
g\,= & \,{2 \pi B \over  \text{sin} \, (2 \pi B z) }\,e^{z(D x + E)}\,,
\tag 4.24
\\
g\,=\,& {2 \pi B \, \text{sin} \, (2 \pi B (x - x_0 - z ))
\over  \text{sin} \, (2 \pi B (x - x_0)) \, \text{sin} \, (2 \pi B z) }
\,e^{z(D x + E)}\,
\tag 4.25
\endalign
$$
where $ B, D, E, x_0$ are parameters.

Elliptic type.
$$
g\,=\, 
B\, 
{\theta (B (x - x_0 -z ), \tau) \, \theta ' (0, \tau) 
\over
\theta (B (x - x_0 ), \tau)\, \theta  ( B z , \tau) }
\, e^{z(D x + E)}
\tag 4.26
$$
where $ B, D, E, x_0, \tau$ are parameters.
\endproclaim

\demo{Proof}
Let $g$ be a limit of $f$. Introduce $v_1 = \partial _z \,
(\partial _x + \partial _z)\, \text{ln}\, g$,
\, $v_2 = \partial _x \,
(\partial _x + \partial _z)\, \text{ln}\, g$,
\, $v_3 =  \partial _x \,
 \partial _z\, \text{ln}\, g$. 

\proclaim{Lemma 4.15} 
The functions $v_1,v_2,v_3$ have the following properties.

1. $v_1$
is a function of $z$, $v_2$
is a function of $x$, $v_3$
is a function of $x-z$. 

2. After identifying the respective variables
in these functions with a new variable $t$, we have
$v_2(t)=v_3(t)$. 

3. $v_1 (t) = t^{-2} + O(t^{-1}),$ if $ t\to 0$.

4. The 
functions $v_1, v_2, v_3$ satisfy a common 
differential equation of the form
$$
(v')^2\,=\, 4 v^3 + p v^2 + q v + r\tag 4.27
$$
 for suitable numbers $p, q, r$. Such an equation is unique. 
\endproclaim

\demo{Proof}
Properties 1,2,3 are satisfied for $g=f$ with any values of the 
parameters. Therefore, they are satisfied for any limit $g$.

 Property 4
is satisfied for $g=f$. Using property 3, 
we conclude that for any limit $g$ 
the function $v_1$ satisfies a unique equation of form (4.27).
Since for $g=f$ the functions $v_2,v_3$ satisfy the same equation as $v_1$, 
this is also true for any limit. The Lemma is proved. $\square$
 \enddemo
 
Let $P(t)$ be the cubic polynomial on 
the right hand side of (4.27).

\proclaim{Lemma 4.16}
Let the roots of $P$ be pairwise distinct. 
Then $g$ has  form (4.26).
\endproclaim

\demo{Proof} If the roots are distinct, then by Lemma 4.15 the function
$v_1(t)$ has the form $B^2\wp(Bt, \tau) + D$
where $\wp$ is the Weierstrass function, 
and $B, D,\tau$ are 
suitable constants. The functions $v_2=v_3$ satisfy the 
same differential equation (4.27) as $v_1$. Therefore, 
either $v_2(t)=v_3(t)=v_1(t-t_0)$ for some $t_0\in\C$, or
$v_2=v_3=\text{const}$. It is clear that the second situation cannot
arise in a limit of $f$, so 
$v_2(t)=v_3(t)=v_1(t-t_0)$.

If $v_1, v_2, v_3 $ are known, then
the second differential of the logarithm
of the function $g$ is known, and hence the function $g$ is known up to a 
transformation of the form $g \mapsto g\,e^{ax + bz + c}$
for suitable constants $a, b, c$. Therefore, 
$g = G \, e^{ax + bz + c}$ where $G$ has  form (4.26). The condition
Res$_{z=0}g = 1$ implies $a=c=0$. Now
$g = G \, e^{ bz}$ and the parameter $b$ can be included into the parameter
$E$ of (4.26). The Lemma is proved. $\square $

\enddemo

\proclaim{Lemma 4.17}
Let $P$ have a root of multiplicity $2$.
Then $g$ has  form (4.24) or (4.25).
\endproclaim

\demo{Proof} An equation of the form 
$(v')^2= 4(v-\al)^2 (v-\beta)$ 
can be solved explicitly.
This gives
$$
v_1 (z)\,=\, {4\pi^2 B^2 \over
\text{sin}^2\, ( 2\pi B z)} \, +\, D\,.
\tag 4.28
$$
The function $v_2=v_3$ has to be a solution of this equation, 
so either $v_2(t)=v_3(t)=v_1(t-t_0)$, or $v_2=v_3=\alpha$, or $v_2=v_3=\beta$. 
It is easy to see that the third case cannot arise as a limiting case of $f$. 
In the first case, 
$g$ has  form (4.25) up to a factor  $e^{ax+bz+c}$
for suitable numbers $a, b, c$. Reasoning as before we conclude that
$g$ has form (4.25).

Similarly, in the second case, 
$g$ has  form (4.24). The Lemma is proved. $\square $

\enddemo

\proclaim{Lemma 4.18} 
If the polynomial $P$ has a root of multiplicity $3$, then
 $g$ has  form (4.22) or (4.23).
\endproclaim

\demo{Proof} Analogous to Lemma 4.17.  $\square $ \enddemo

The Proposition is proved. $\square $

\enddemo

{\bf 4.5.5. End of proof of Theorem 4.2, part 2.} 

In the previous section 
we have determined the possible forms of the function $\mu(x,z)$
for any root $\alpha$: they are given by (4.22)-(4.26) 
for coupling constant 1. Now we will determine the consistency conditions,
which are imposed on these functions for different roots 
by the $\g_\al\o\g_\beta\o\g_\gamma$-part of the CDYB equation, where
$\al+\beta+\gamma=0$. 

First of all, 
by our assumptions the function $\phi_\al$ is a meromophic function
for any root $\al$. Since $\phi_\al=\mu e^{zL_\al \psi}$, 
and $\mu$ is meromorphic, 
the function $\psi (\la)$ is holomorphic on $\h^*$.
Therefore, by using gauge transformations of type 1 and 2, it is possible 
to reduce any r-matrix with spectral parameter
and coupling constant 1 to the form in which its $\h\o\h$-part
is $T(z)\sum x_i\o x_i$, where $T(z)=\frac{1}{z}+T_1z+O(z^3)$ is an odd,
scalr-valued meromorphic 
function. We will call such $r$-matrices reduced, and from now
on will work only with them. 

For a reduced r-matrix  
$\phi_\alpha (\la, z) = \mu_\al^* ((\al, \la), z)$, where
$\mu_\al^*(x,z)=\mu_\al(x,z)e^{T_1xz}$,
and $\mu_\al$ is the function $\mu$ introduced in Section 4.5.2. 
Observe that $\mu_\al^*$ is a function from family (4.21).
Let $A_\al,B_\al,C_\al,D_\al,E_\al,\tau_\al$ be the parameters
$A,B,C,D,E,\tau$ determined from (4.21), for $\mu=\mu_\al^*$. 
Since the coupling constant is 1, we have $A=-B$. 

\proclaim{Lemma 4.19} All $\mu_\al^*$ are of the same type 
(rational, trigonometric, or elliptic). 
\endproclaim

\demo{Proof}
 From equation (4.16) we can find the function $t(u)$. 
We know that it is the same for all roots $\alpha$. It is easy to check that 
$\mu_\alpha^*$ is of rational, trigonometric, or elliptic type
 iff the set of poles of $t(z)$ is a lattice of rank $0$, $1$, $2$ 
respectively. 
The Lemma is proved. $\square$
\enddemo

So it remains to consider the rational, trigonometric, and elliptic cases 
separately. 

{\bf Elliptic case.}

\proclaim{ Lemma 4.20} Let $r$ be a reduced dynamical r-matrix
such that the functions $\mu^*_\alpha$ are of elliptic type. 
Then there exist complex numbers $a, b, \tau$, Im $\tau>0$,
 and elements
$\nu, \kappa \in \h^*$ such that
$$
A_\al=a, \qquad
B_\al=b, \qquad
\tau_\al=\tau, \qquad
C_\al=(\al,\nu), \qquad
D_\al=c,\qquad
E_\al=(\al,\kappa),
$$
for all roots $\al$. 

\endproclaim
\demo{Proof}
Substitute formula (4.21) into 
the $\g_\al \otimes \g_{\bt} \otimes \g_{\gm}-$part of the CDYB equation,
see (4.17).
Considering the poles at the hyperplanes 
of the form $z_{1,2} =$const, we
observe that the functions $\phi_\bt$ and $\phi_\gm$
have the same lattice of periods. This allows us to conclude that
there exist numbers $b$ and $\tau$ such that
$B_\al=b,$ and $\tau_\al=\tau$ for all $\al$.
The residue condition (4.3) implies the existence of a number $a$ such
that $A_\al=a$ for all $\al$. 
The necessity to cancel the poles at the hyperplanes of the form
$(\al, \la)=$const implies the existence of elements
$\kappa, \nu \in \h^*$, $c\in \Bbb C$  such that
$C_\al=(\al,\nu), D_\al=c, E_\al=(\al,\kappa)$ for all $\al$.

This proves the Lemma, and Theorem 4.2, part 2, in the elliptic case. 
$\square$
\enddemo

Now consider the trigonometric and rational case. 

\proclaim{Lemma 4.21} $D_\alpha$ are the same for all $\alpha$. 
\endproclaim

\demo{Proof} Easily follows from (4.17). $\square$
\enddemo

Thus, we can reduce to the situation $D_\alpha=0$ by using a gauge 
transformation of type 2 with $\psi=-D(\la,\la)/2$.

{\bf Rational case.}

\proclaim{Lemma 4.22} In the rational case 
$E_{\alpha+\beta}=E_\alpha+E_\beta$. 
\endproclaim

\demo{Proof} Follows directly from (4.17). $\square$
\enddemo

Thus, in the rational case we can reduce to the situation
$E_\al=0$ by a gauge trasformation of type 2 with $\psi=-(E,\lambda)$,
$E\in \h^*$. 

If $r$ is a reduced rational dynamical 
r-matrix with $D_\alpha=E_\alpha=0$, then it 
is easy to see 
from Proposition 4.14 that $r$ is of the form $\frac{\Omega}{z}+r_0(\la)$, 
where $r_0$ is skew-symmetric. Since $r_0=\lim_{z\to\infty}r$, 
$r_0(\la)$ is a classical dynamical r-matrix 
without spectral parameter with zero coupling constant. 
Such r-matrices were classified in Theorem 3.2. Theorem 3.2. implies 
that $r(\la,z)$ is equivalent to (4.9) by gauge transformations 
of type 
$3$. 
This proves Theorem 4.2, part 2, in the rational case. 

{\bf Trigonometric case.} 
In the trigonometric case,
the functions $\mu_\alpha^*$ have form (4.24),(4.25). 
As in the elliptic case, 
it is easy to see that $B$ is the same for all $\alpha$
since $B^{-1}$ is the period of the lattice of poles. So by a gauge 
transformation of type 4 we can arrange $B=1$.     

\proclaim{Lemma 4.23} Let $\rho$ be the half-sum of positive roots. 
There exists a limit of $r(s\rho,z)$ as $s\to i\infty$
\endproclaim

\demo{Proof} 
The statement is clear from formulae (4.24),(4.25).
\enddemo

Denote by $\bar r(z)$ this limit. This is a classical r-matrix  
with  spectral parameter having the form $\frac{\Omega}{z}+O(1)$
as $z \to 0$ and 
invariant under the action of the
 Cartan subalgebra. A classification of such r-matrices
was given by  Belavin and Drinfeld  \cite{BD1}, Theorem 6.1.

Namely, consider an  r-matrix 
$$
r_{tr}(z)=2i\frac{\Omega_-e^{2iz}+\Omega_+}{e^{2iz}-1},
\tag 4.29
$$
where 
$\Omega_\pm=\frac{1}{2}\sum x_i\o x_i +
\sum_{\alpha\in\Delta_\pm}e_\alpha\o e_{-\alpha}$ are the half-Casimirs.
According to \cite{BD1}, any classical r-matrix of the above type 
 can be obtained from $r_{tr}(z)$ 
by change of polarization of the Lie algebra, and by gauge transformations 
of type 
$4$, and type $2$ with a linear function $\psi$.

Thus, in order to prove Theorem 4.2, part 2, in the trigonometric case, 
it is enough to assume that $\bar r(z)=r_{tr}(z)$. 

Under this assumption, it is easy to deduce from Proposition 4.14 
that $\mu_\alpha^*(x,z)$ equals to:

1. $\frac{\sin(x-x_0+z)}{\sin(x-x_0)\sin z}$
if $\mu_\al^*$ depends nontrivially on $x$; 

2. $\frac{e^{-iz}}{\sin z}$, if $\al>0$ and $\mu_\al^*$ does not depend on
$x$. 

3. $\frac{e^{iz}}{\sin z}$, if $\al<0$ and $\mu_\al^*$ does not depend on
$x$.

Now let us send $z\to i\infty$. It is easy to see that the limit of 
$r(\la,z)$ exists. Denote this limit by $r_\infty(\la)$. Let 
$\mu_\al^\infty((\la,\al))$ be the coefficient of $e_\alpha\o e_{-\alpha}$ 
in $r_\infty (\la)$. From cases 1-3 above, we get
that $\mu_\al^\infty(x)$ equals to:

1'. $\frac{e^{-i(x-x_0)}}{\sin(x-x_0)}$,
if $\mu_\al^\infty$ depends nontrivially on $x$; 

2'. $-2i$, if $\al>0$ and $\mu_\al^\infty$ 
does not depend on
$x$. 

3'. $0$, if $\al<0$ and $\mu_\al^\infty$ does not depend on
$x$.

On the other hand, $r_\infty(\la)$ is 
a classical dynamical r-matrix without  spectral parameter, and 
it is clear from cases 1'-3' that it has coupling constant $-2i$.  

So, $r_\infty(\la)$ has to be of the form given by Theorem 3.10.
This determines possible 
combinations of functions $\mu_\alpha^*$, and shows that 
$r(\la,z)$ is equivalent to an r-matrix (4.8) for a suitable subset $X$
by a gauge transformation of type 3. 

Theorem 4.2, part 2 is proved. 

\subhead 4.6. Dynamical r-matrices without spectral parameter 
for affine Lie algebras
\endsubhead

In this section we  interpret dynamical classical r-matrices 
without spectral parameter for an affine Lie algebra as 
dynamical r-matrices with spectral parameter for the underlying
simple Lie algebra. 

Let $\g$ be a simple Lie algebra, and $\tilde\g=\g[t,t^{-1}]\oplus 
\C c\oplus\C d$ be the corresponding affine Lie algebra,
where $c$ is the central element, and $d$ is the grading element. 
Let $\h \subset \g$ be a Cartan subalgebra, and $\{x_i\}$ an orthonormal basis 
of $\h$. Let $\tilde\h\subset\tilde\g$ be the Cartan subalgebra
of $\tilde\g$, $\tilde\h=\h\oplus \C c\oplus \C d$.
Recall that $c,d$ are orthogonal to $\h$ with respect to the 
standard bilinear form, and $(c, d) = 1$, $(c, c) = (d, d) = 0$. 
 
The elements of $\tilde \h$ have the form $h + xc + yd$ where $h \in \h, \,
x, y \,\in \, \C$. The elements of $\tilde \h^*$ are triples $(\la, k, s )$
such that $\< (\la, k, s ), h + xc + yd \>\,=\, \<\la , h \> + kx + sy.$

Let $\dl \in \tilde \h^*$ be the positive imaginary root, $\< \dl, d\> = 1,
\< \dl, c\> = 0, \< \dl, \h \> = 0.$  
The roots of $\tilde\g$ are $\alpha+n\delta$ and $n \dl$, where
$\alpha$ is a root of $\g$, 
$n \in \Z$, and in the second case $n \neq 0$. The positive roots are
$\alpha+n\delta$ and $(n+1) \dl$, where $\al$ is a positive root of $\g$, and
 $n \geq 0$.

A basis of positive root elements of $\tilde\g$ is formed by
$e_\al t^n$, $e_{-\al}t^{n+1}$, and $x_it^{n+1}$, where $\al>0, \, n  \ge 0$. 
The dual elements are respectively
$e_{-\al} t^{-n}$, $e_{\al}t^{-n-1}$, and $x_it^{-n-1}$, where $\al>0, n\ge 0$. 
  
According to Theorem 3.1, the solution  
of CDYBE for 
$\tilde\g$ with coupling constant $2$ and $C=0,\nu=0$ has the form
$$
r(\la,k,s)=\hat \Omega+\sum_{\al\in\Delta, \, n}
\text{cotanh}((\al,\la)+kn)e_\al t^n\o e_{-\al}t^{-n}+
\sum_{i}\sum_{n\ne 0}\text{cotanh}(kn)x_it^n\o x_it^{-n}
.\tag 4.23
$$
where  $\hat\Omega$ is the  Casimir element of $\tilde \g$. 

For any $z\in \C^*$,  let $\pi_z: \g[t,t^{-1}]\to \g$ be the evaluation 
map at $t=z$. Consider 
the function $\bar r(y,w)=(\pi_y\o \pi_w)(r|_{c=0})$,
where $r|_{c=0}\in \g[t,t^{-1}]^{\hat\o 2}$ is the image of $r$ 
under the reduction modulo $c$. It is easy to see that $\bar r$ 
depends only on $u=y/w$, and equals
$$
\gather
\bar r(\la,k,u)=\\
\sum_{\al\in\Delta}[\sum_{n\in\Bbb Z}
u^n(1+\text{cotanh}((\al,\la)+kn))]e_\al\o e_{-\al}+
\sum_{i}[1+\sum_{n\ne 0}u^n(1+\text{cotanh}(kn))]x_i\o x_i
,\tag 4.24
\endgather
$$

Set $u=e^{2\pi iz}$.

\proclaim{Lemma 4.24} For any $k$, the function 
$\hat r(\la,z)=\bar r(\la,k,u)$ satisfies the CDYB equation with  spectral 
parameter $z$. 
\endproclaim

\demo{Proof} Applying the operator $\pi_{z_1}\o \pi_{z_2}\o \pi_{z_3}$ 
to the CDYB equation (3.3) for $\tilde \g$, one easily
obtains the CDYB equation with spectral parameter. $\square$
\enddemo

Now we 
compute $\hat r(\la,z)$. Set $\tau=k/\pi i$, assume that Im $\tau > 0$,
and use  the classical formulae
$$
\sum_{n\in \Bbb Z}u^n(1+ \text{cotanh}(a+\pi i\tau n))=-\frac{1}{\pi i}
\sigma_{\frac{a}{\pi i}}(z,\tau),\tag 4.30
$$
and
$$
1+\sum_{n\in \Bbb Z\setminus 0}u^n(1+ \text{cotanh}(\pi i\tau n))=-\frac{1}{\pi i}
\rho(z,\tau),
\tag 4.31
$$
where $\sigma$ and $\rho$ are defined in
(4.6). Then  $\hat r(\la,z)$ takes the form of the Felder
solution of CDYBE with spectral parameter,
$$
\hat r(\la, z)\,=\, -\, \frac{1}{\pi i}\,\rho(z,\tau)\, \sum_{i=1}^N
x_i \otimes x_i\,-\, \frac{1}{\pi i}\,
\sum_{\al \in \Dl}\,\sigma_{\frac{(\al, \la) }{\pi i}}(z,\tau)
\, e_\al \otimes e_{-\al}\,.
$$

\subhead{ Appendix: Open problems}\endsubhead

In conclusion we would like to formulate two open problems. 

Let $\frak g$ be a simple Lie algebra, $\frak h$ its Cartan subalgebra, 
and $\h_0$ a subspace in $\frak h$. 
We will say that a meromorphic function $r:\h_0^*\to \g\o\g$ 
is a classical dynamical r-matrix if it satisfies (3.2),(3.3), and (3.1) for 
$h\in\h_0$. We will say that a meromorphic function
$r:\h_0^*\times \C\to \g\o\g$ is a classical dynamical r-matrix 
with a spectral parameter if it satisfies (4.2)-(4.4), and 
(4.1) for $h\in \h_0$. 
The number $\epe$ in both cases is called the coupling constant. 

\proclaim{Problem 1} Classify classical dynamical r-matrices on $\h_0$, 
with a nonzero coupling constant. 
\endproclaim

\proclaim{Problem 2} Classify classical dynamical r-matrices on $\h_0$ 
with a spectral parameter, 
with a nonzero coupling constant. 
\endproclaim

These problems are solved for two extreme cases -- $\h_0=0$ (\cite{BD1,BD2})
and $\h_0=\h$ (this paper). We expect that it can be solved 
in the intermediate cases by combining of the methods of our paper 
and the two papers of Belavin and Drinfeld. 

{\bf Note added in the second version:} While this paper was
being revised, O.Schiffmann obtained a partial solution 
of Problem 1 \cite{Sch}.

\vskip10ex 

P.Etingof
Department of Mathematics, Harvard University, Cambridge, MA 02138

{\it E-mail address: etingof\@math.harvard.edu.}

A.Varchenko,
Department of Mathematics, Phillips Hall, University of North Carolina at
Chapel Hill, Chapel Hill, NC 27599-3250, USA

{\it E-mail address: av\@math.unc.edu.}

\end